\documentclass{amsart}
\usepackage{amsmath,amssymb}
\usepackage[dvips]{graphics}
\usepackage[all]{xy}
\newtheorem{Theorem}{Theorem}[section]
\newtheorem{Lemma}[Theorem]{Lemma}

\newtheorem{Corollary}[Theorem]{Corollary}

\theoremstyle{Definition}

\theoremstyle{Remark}
\newtheorem{Remark}[Theorem]{Remark}

\makeatletter
\@addtoreset{figure}{section}%{subsection}
\def\@thmcountersep{-}
\makeatother

\numberwithin{equation}{section}

\begin{document}

\title{Regular projections of graphs with at most three double points}

\author{Youngsik Huh}
\address{Department of Mathematics, School of Natural Sciences, Hanyang University, Seoul 133-791, Korea}
%\curraddr{}
\email{yshuh@hanyang.ac.kr}
\thanks{The first author was supported by the Korea Science and Engineering Foundation (KOSEF) grant funded by the Korea government (MOST) (No.~R01-2007-000-20293-0).}

\author{Ryo Nikkuni}
\address{Department of Mathematical Sciences, School of Arts and Sciences, Tokyo Woman's Christian University, 2-6-1 Zempukuji, Suginami-ku, Tokyo 167-8585, Japan}
\email{nick@lab.twcu.ac.jp}
\thanks{The second author was partially supported by Grant-in-Aid for Young Scientists (B) (No. 18740030), Japan Society for the Promotion of Science.}

%    General info

\date{}

\dedicatory{}


\begin{abstract}
%We investigate some relations between a generic immersion of a graph into the 2-space with small number of double points and the embeddings of the graph into the 3-space obtained by lifting the immersion with respect to the natural projection from the 3-space to the 2-space. Specifically we prove:
%(1) An embedding of a graph obtained from a generic immersion of the graph with at most three double points is totally free if it contains neither a Hopf link nor a trefoil knot.
%(2) If a generic immersion of a planar graph does not produce any trivial embedding of the graph, then the number of double points of the immersion is more than or equal to three.
%
A generic immersion of a planar graph into the 2-space is said to be knotted if there does not exist a trivial embedding of the graph into the 3-space obtained by lifting the immersion with respect to the natural projection from the 3-space to the 2-space. In this paper we show that if a generic immersion of a planar graph is knotted then the number of double points of the immersion is more than or equal to three. To prove this, we also show that an embedding of a graph obtained from a generic immersion of the graph (does not need to be planar) with at most three double points is totally free if it contains neither a Hopf link nor a trefoil knot.
\end{abstract}

\maketitle

\section{Introduction}

Throughout this paper we work in the piecewise linear category and graphs are considered as topological spaces. Let ${\mathbb S}^{3}$ be the unit $3$-sphere in ${\mathbb R}^{4}$ centered at the origin. For a finite graph $G$, an embedding $f:G \to {\mathbb S}^{3}$ is called  a {\it spatial embedding of $G$} or simply a {\it spatial graph}. If $G$ is homeomorphic to the disjoint union of $n$ circles, then $f$ is called an {\it $n$-component link} (or a {\it knot} if $n=1$). Two spatial embeddings $f$ and $g$ of $G$ are said to be {\it equivalent} ($f \approx g$) if there exists an orientation-preserving self-homeomorphism $\Phi$ on ${\mathbb S^{3}}$ such that $\Phi(f(G))=g(G)$. A graph $G$ is said to be {\it planar} if there exists an embedding of $G$ into the unit $2$-sphere ${\mathbb S}^{2}$. A spatial embedding of a planar graph $G$ is said to be {\it trivial} if it is equivalent to an embedding $h:G \to {\mathbb S}^{2}\subset {\mathbb S}^{3}$.

\begin{figure}[htbp]
\[
\xymatrix{
\scalebox{0.35}{\includegraphics*{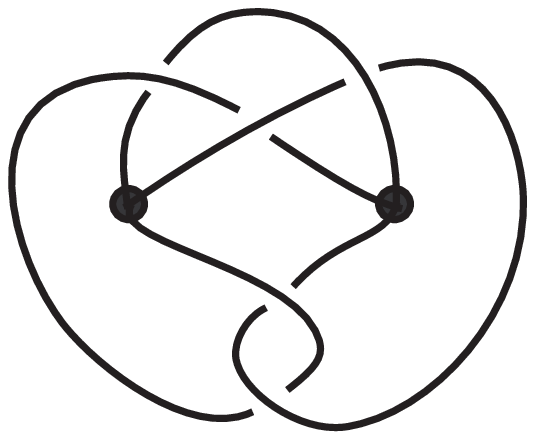}}
\ar[r]^{\approx} &
\scalebox{0.35}{\includegraphics*{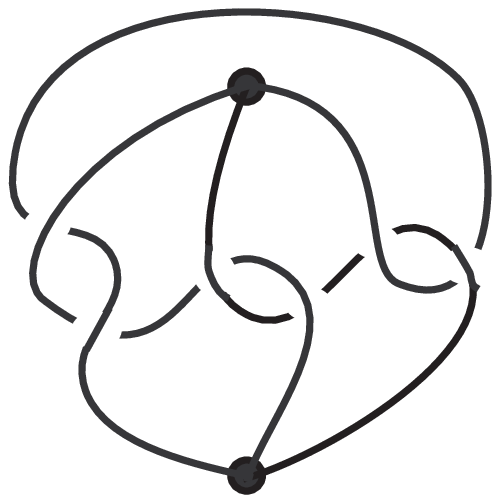}}
\ar[d]^{\pi} \\
\scalebox{0.35}{\includegraphics*{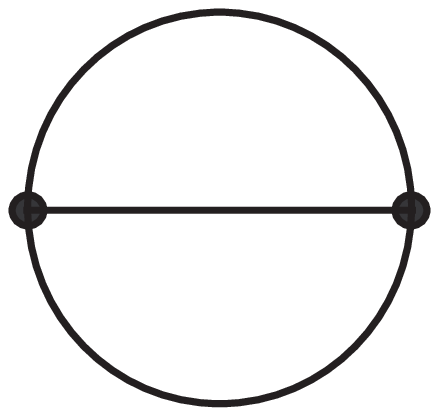}}
\ar[u]^{f}  \ar[ru]^{f'} \ar[r]_{\varphi} &
\scalebox{0.35}{\includegraphics*{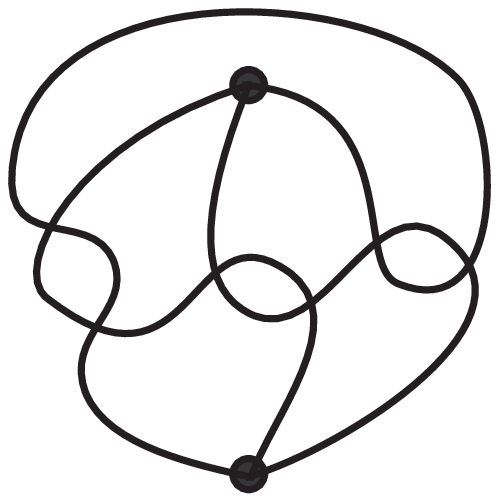}}
}
\]
   \caption{}
  \label{projects}
\end{figure}

A continuous map $\varphi:G \to {\mathbb S}^{2}$ is called a {\it regular projection} of $G$ if the multiple points of $\varphi$ are only finitely many transversal double points away from the vertices of $G$. For a spatial embedding $f$ of $G$, we also say that $\varphi$ is a {\it regular projection of $f$} or {\it $f$ projects on $\varphi$}, if there exists an embedding $f':G \to {\mathbb S}^{3}\setminus \left\{(0,0,0,1),(0,0,0,-1)\right\}$ such that $f$ is equivalent to $f'$ and $\pi\circ f'=\varphi$, where $\pi:{\mathbb S}^{3}\setminus \left\{(0,0,0,1),(0,0,0,-1)\right\}\to {\mathbb S}^{2}$ is the natural projection, see Fig. \ref{projects}. A {\it regular diagram $\widetilde{\varphi}$ of $f$} is none other than the regular projection $\varphi$ of $f$ with over/under information of  each double point. We call a double point with over/under information a {\it crossing}. For a subspace $H$ of $G$, we often denote $\varphi(H)$ (resp. $\widetilde{\varphi}(H)$) by $\widehat{H}$ (resp. $\widetilde{H}$) as long as no confusion occurs.

\begin{figure}[htbp]
\[
\xymatrix{
\scalebox{0.4}{\includegraphics*{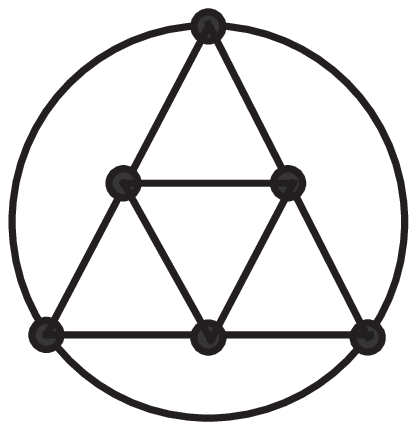}}
\ar[r]^{\varphi} &
\scalebox{0.4}{\includegraphics*{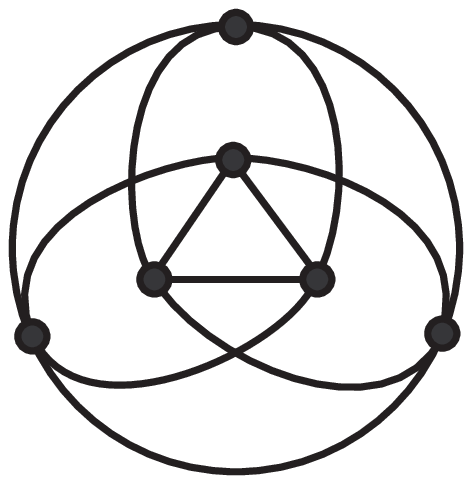}}
}
\]
   \caption{}
  \label{kp}
\end{figure}

Our purpose in this paper is to investigate knotted projections which realize the minimal number of double points.
A regular projection $\varphi$ of a planar graph $G$ is said to be {\it knotted} if there does not exist any trivial spatial embedding of $G$ which projects on $\varphi$. Such a regular projection was  discovered by K. Taniyama first \cite{Tani95}. For example, let $\varphi$ be the regular projection illustrated in Fig. \ref{kp}. Then we can see that any spatial embedding of $G$ which projects on $\varphi$ contains a Hopf link, so $\varphi$ is knotted. We call a knotted regular projection simply a {\it knotted projection}. By the notion of knotted projection, a problem in graph minor theory can be formulated \cite{Lovasz06}. A planar graph is said to be {\em trivializable} if it has no knotted projections. Let $\Omega$ be the set of all non-trivializable planar graphs whose all proper minors are trivializable. It is known that for any trivializable planar graph $G$ every minor of $G$ is also trivializable \cite{Tani95}. Therefore, due to the celebrated work of Robertson and Seymour on graph minors \cite{RS}, it is guaranteed that $\Omega$ is finite. But, although many elements of $\Omega$ have been found out through continued works \cite{SS00, Tamu04, NOTT05, N07}, the set is not completely determined yet. In this paper, as an effort on this issue,  we will give necessary conditions for knotted projections. 

Let $d$ be a double point of a regular projection $\varphi$ of $G$ such that $\varphi^{-1}(d)=\left\{p_{1},p_{2}\right\}$ and $p_{i}\in e_{i}$, where $e_{i}$ is an edge of $G$ ($i=1,2$). Then we say that $d$ is {\it Type-S} if $e_{1}=e_{2}$, {\it Type-A} if $e_{1}\neq e_{2}$ and $e_{1}\cap e_{2}\neq \emptyset$, and {\it Type-D} if $e_{1}\cap e_{2}=\emptyset$. Then we have the following. Here we denote the number of all double points of a regular projection $\varphi$ by ${\rm cr}(\varphi)$. And a regular projection of a planar graph $G$ is said to be {\it trivial} if only trivial spatial embeddings of $G$ project on it

\begin{Theorem}\label{main2}
Let $\varphi$ be a regular projection of a planar graph $G$. Then we have the following.
\begin{enumerate}
\item If ${\rm cr}(\varphi)=1$, then $\varphi$ is trivial.
\item If ${\rm cr}(\varphi)=2$, then $\varphi$ is not knotted. Moreover, $\varphi$ is trivial if $\varphi$ has a double point of Type-S or Type-A.
\item If ${\rm cr}(\varphi)=3$, then $\varphi$ is not knotted if $\varphi$ has a double point of Type-S or Type-A.
\end{enumerate}
\end{Theorem}

As a corollary of Theorem \ref{main2}, necessary conditions for knotted projections are derived.

\begin{Corollary} \label{main2cor}
If a regular projection $\varphi$ of a planar graph is knotted, then ${\rm cr}(\varphi)\ge 3$. In partcular, if $\varphi$ is knotted and ${\rm cr}(\varphi)=3$ then every double point of $\varphi$ is Type-D.
\end{Corollary}

As we saw in Fig. \ref{kp}, there exists a knotted projection with only three double points. Thus the inequality of Corollary \ref{main2cor} is best possible. 

\begin{figure}[htbp]
      \begin{center}
\scalebox{0.35}{\includegraphics*{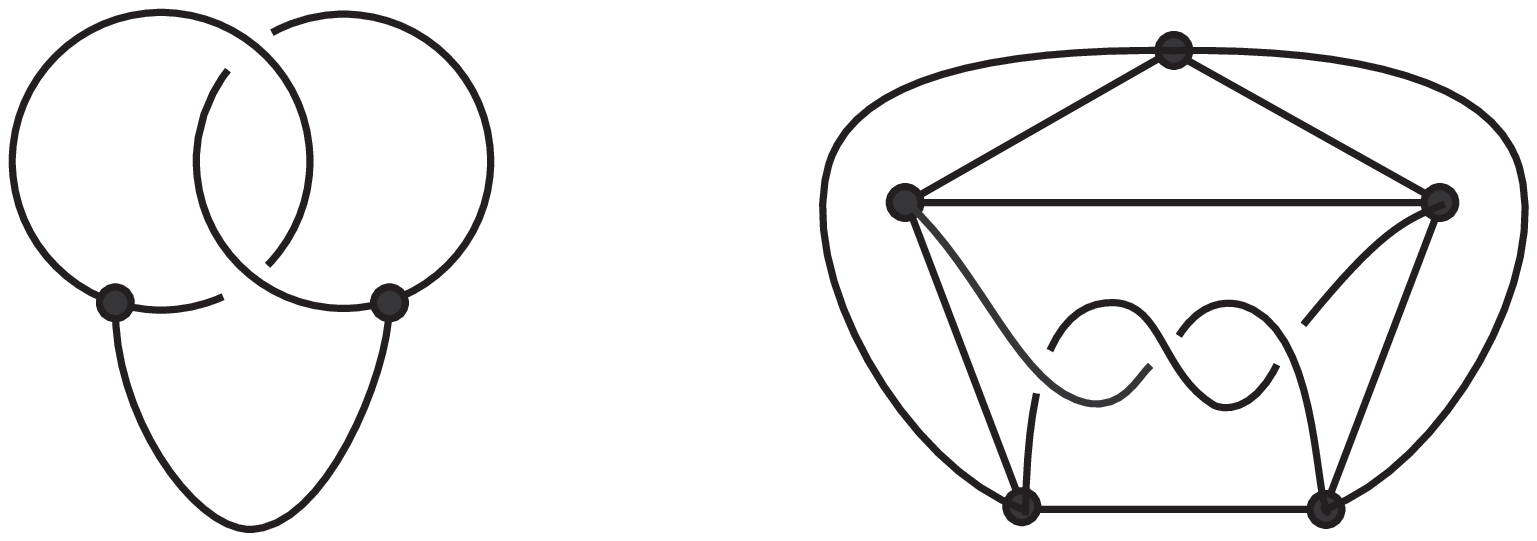}}
      \end{center}
   \caption{}
  \label{not_totally_free}
\end{figure}

To accomplish the proof of Theorem \ref{main2}, we determine non-trivial spatial graph types which may be contained in a spatial embedding of a graph which projects on a regular projection of the graph with at most three double points. In the case of an $n$-component link $L$ which projects on a regular projection with at most three double points, it is not hard to see in knot theory that $L$ is trivial if it does not contain a Hopf link or a trefoil knot. In the following we generalize the above fact to spatial graphs. 
A spatial embedding $f$ of a graph $G$ is said to be {\it free} if the fundamental group of the spatial graph complement $\pi_{1}({\mathbb S}^{3}\setminus f(G))$ is free. Moreover we say that $f$ is {\it totally free} if the restriction map $f|_{H}$ is free for any subgraph $H$ of $G$. For example, the two spatial graphs in Fig. \ref{not_totally_free} are free but not totally free. We remark here that if $G$ is planar, then $f$ is totally free if and only if $f$ is trivial by  Scharlemann-Thompson's famous theorem \cite{ST91}. Then we have the following. 

\begin{Theorem}\label{main1}
Let $\varphi$ be a regular projection of a graph $G$ and $f$ a spatial embedding of $G$ which projects on $\varphi$. Assume that ${\rm cr}(\varphi)\le 3$. Then $f$ is totally free if it does not contain a Hopf link or a trefoil knot.
\end{Theorem}

As a direct consequense of Theorem \ref{main1}, we have the following corollary.

\begin{Corollary}\label{main1cor}
Let $\varphi$ be a regular projection of a planar graph $G$ and $f$ a spatial embedding of $G$ which projects on $\varphi$. Then we have the following.
\begin{enumerate}
\item If ${\rm cr}(\varphi)=1$, then $f$ is trivial.
\item If ${\rm cr}(\varphi)=2$, then $f$ is trivial if it does not contain a Hopf link.
\item If ${\rm cr}(\varphi)=3$, then $f$ is trivial if it does not contain a Hopf link or a trefoil knot.
\end{enumerate}
\end{Corollary}

Corollary \ref{main1cor} also leads to another fundamental result on spatial graphs. A spatial embedding $f$ of a planar graph $G$ is said to be {\it minimally knotted} if $f$ is not trivial but $f|_{H}$ is trivial for any proper subgraph $H$ of $G$.  Fig. \ref{projects} shows an example of minimally knotted spatial embedding which is called {\it Kinoshita's theta curve}. Note that every planar graph without isolated vertices and free vertices has minimally knotted spatial embeddings \cite{Kawauchi89, Wu93}.

\begin{Corollary}\label{main1cor2}
Let $\varphi$ be a regular projection of a planar graph $G$ and $f$ a minimally knotted spatial embedding of $G$ which projects on $\varphi$. If $f$ is neither a Hopf link nor a trefoil knot, then ${\rm cr}(\varphi)\ge 4$.
\end{Corollary}

\begin{proof}
If ${\rm cr}(\varphi)\le 3$, by Theorem \ref{main1} we have that $f$ contains a Hopf link or a trefoil knot. Since $f$ is minimally knotted, $f$ must be a Hopf link or a trefoil knot.
\end{proof}

\begin{figure}[htbp]
      \begin{center}
\scalebox{0.4}{\includegraphics*{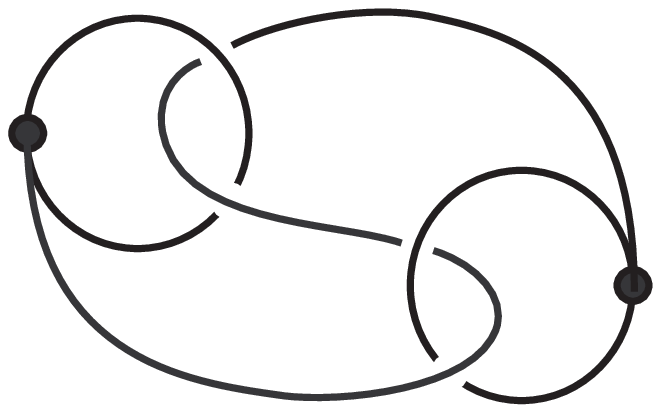}}
      \end{center}
   \caption{}
  \label{handcuff}
\end{figure}

The inequality of Corollary \ref{main1cor2} is best possible. For example, the spatial handcuff graph in Fig. \ref{handcuff} is minimally knotted and it can project on a regular projection with four double points.

By Corollary \ref{main1cor2}, we also give a partial answer for Ozawa's question \cite[Question 3.7]{NOTT05} which asks whether a minimally knotted spatial embedding of a planar graph can project on a knotted projection of the graph or not.

\begin{Corollary}\label{cor3}
Let $\varphi$ be a knotted projection of a planar graph $G$ with ${\rm cr}(\varphi)=3$. Then there does not exist a minimally knotted spatial embedding of $G$ which projects on $\varphi$.
\end{Corollary}

\begin{figure}[htbp]
      \begin{center}
\scalebox{0.4}{\includegraphics*{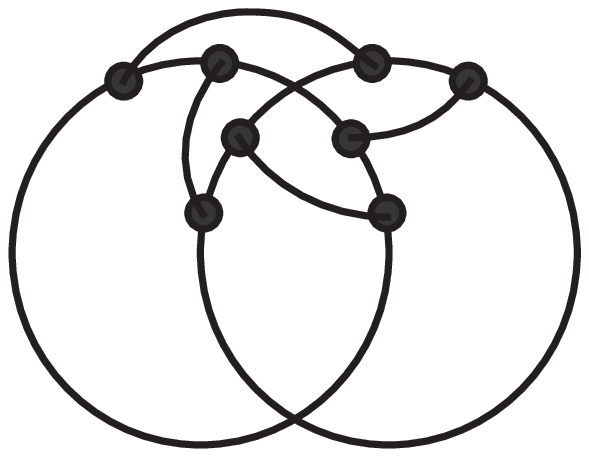}}
      \end{center}
   \caption{}
  \label{TT_ex}
\end{figure}

\begin{Remark}\label{rem}
{\rm 
There exists a regular projection $\varphi$ of a non-planar graph $G$ with ${\rm cr}(\varphi)=2$ such that no totally free spatial embeddings of $G$ project on $\varphi$. For example, let $\varphi$ be the regular projection of a non-planar graph $G$ as illustrated in Fig. \ref{TT_ex}. Then we can see that any of the spatial embedding of $G$ which projects on $\varphi$ contains a Hopf link \cite[Fig. 4]{TT96}. This says that the planarity of a graph is essential in Theorem \ref{main2}. 
}
\end{Remark}

The rest of this paper is organized as follows. In the next section, we introduce a key theorem which is needed to prove our theorems. We prove Theorem \ref{main1} in section 3. And, by utilizing the theorem, the proof of Theorem \ref{main2} is given in section 4. 

\section{Key theorem}

To prove Theorem \ref{main1} and \ref{main2}, we take advantage of a nice geometric characterization of totally free  spatial embeddings which was first proved by Wu for planar graphs \cite[THEOREM 2]{Wu92} and generalized to arbitrary graphs (not need to be planar) by Robertson-Seymour-Thomas \cite[(3.3)]{RST3}. Here a {\it cycle} of a graph $G$ is a subgraph of $G$ which is homeomorphic to the circle, and a {\it disk} is a topological space which is homeomorphic to the unit $2$-disk ${\mathbb D}^{2}$ in ${\mathbb R}^{2}$.

\begin{Theorem}\label{keythm} \cite[(3.3)]{RST3}
A spatial embedding $f$ of a graph $G$ is totally free if and only if for any cycle $\gamma$ of $G$ there exists a disk $D_{\gamma}$ in ${\mathbb S}^{3}$ such that
\begin{eqnarray*}
f(G)\cap D_{\gamma} = f(G) \cap \partial D_{\gamma} = f(\gamma).
\end{eqnarray*}
\end{Theorem}

Namely $f$ is totally free if and only if for any cycle $\gamma$ of $G$ the knot $f(\gamma)$ bounds a disk $D_{\gamma}$ in ${\mathbb S}^{3}$ as a Seifert surface such that ${\rm int}D_{\gamma}\cap f(G)=\emptyset$. We call $D_{\gamma}$ a {\it trivialization disk} for $f(\gamma)$. Theorem \ref{keythm} helps us to detect the totally freedom (or triviality) of a spatial graph by utilizing local informations in the regular diagram.

To put it into practice, we introduce some definitions. Let $\widetilde{\varphi}$ be a regular diagram of a spatial embedding of a graph $G$. Fix a cycle $\gamma$ of $G$. Among the edges of $G$ not contained in $\gamma$, choose all possible edges $e_{1},e_{2},\ldots,e_{m}$ so that $\widehat{\gamma}$ and $\widehat{e}_{i}$ produce double points of $\varphi$. We denote the subgraph of $G$ which is obtained from $G$ by forgetting $e_{1},e_{2},\ldots,e_{m}$ by $G'$. Let $R_{1},R_{2},\ldots,R_{k}$ be all of the connected components of ${\mathbb S}^{2}\setminus \widehat{\gamma}$. We denote the subspace $\varphi^{-1}\left(\widehat{G'}\cap R_{i}\right)$ of $G'$ by $H_{i}\ (i=1,2,\ldots,k)$.

\begin{figure}[htbp]
      \begin{center}
\scalebox{0.35}{\includegraphics*{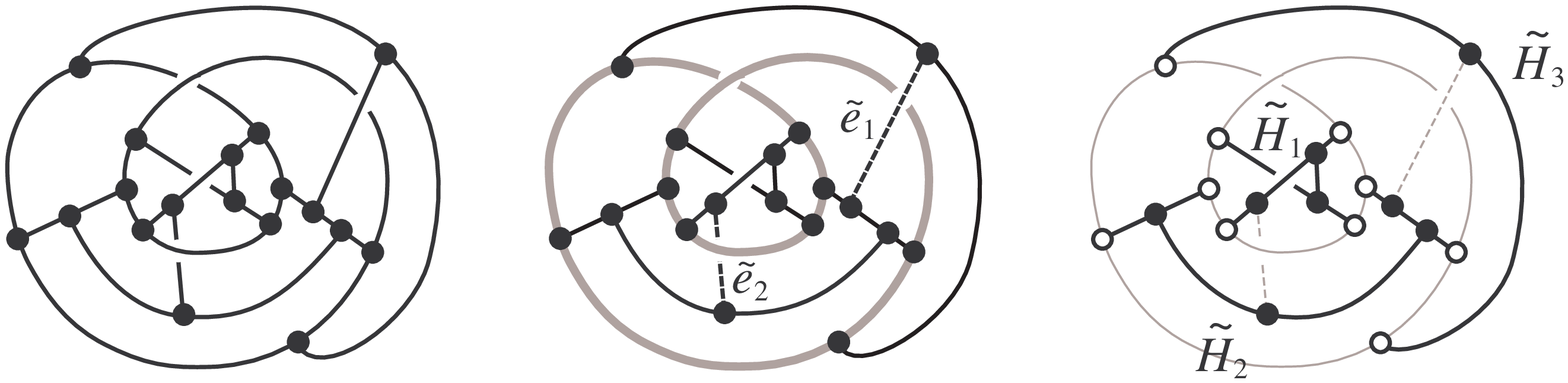}}
      \end{center}
   \caption{}
  \label{local_cycle}
\end{figure}

\begin{figure}[htbp]
      \begin{center}
\scalebox{0.35}{\includegraphics*{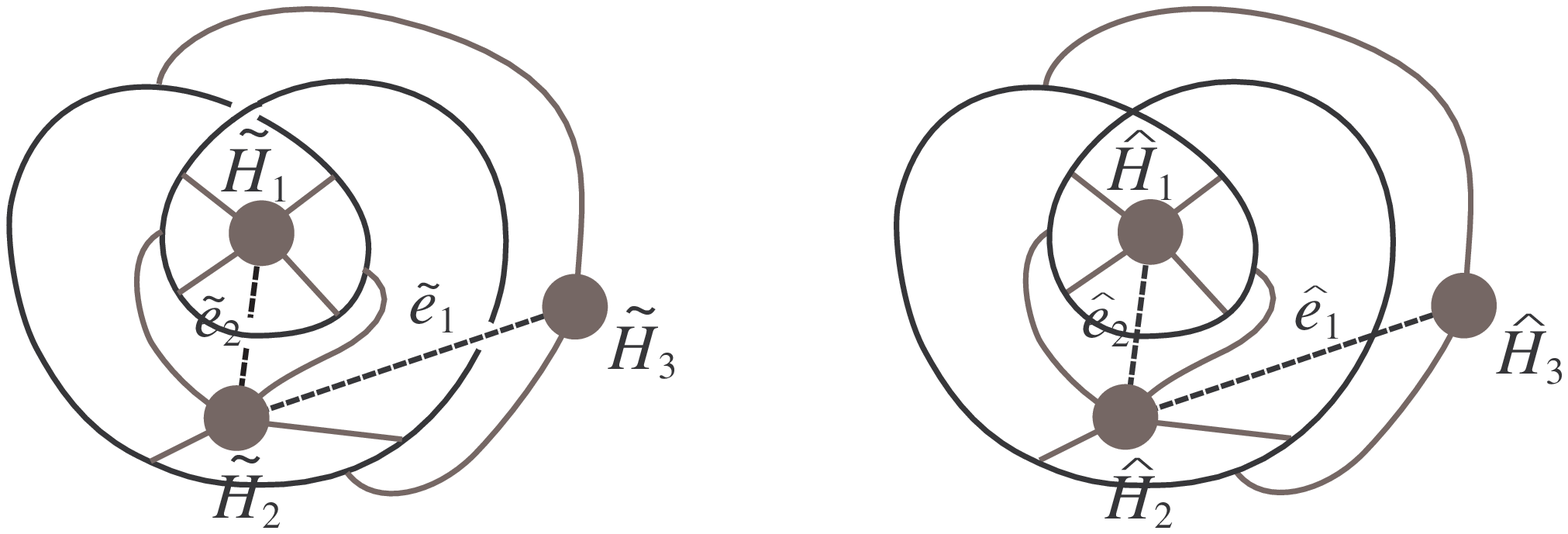}}
      \end{center}
   \caption{}
  \label{local_cycle2}
\end{figure}

For example, given a regular diagram as the left-hand side of Fig. \ref{local_cycle}, let $\gamma$ be the cycle of $G$ such that $\widetilde{\gamma}$ corresponds to the gray curve in the center of Fig. \ref{local_cycle}, where $\widetilde{e}_{1}$ and $\widetilde{e}_{2}$ are drawn by dotted black lines. Then the right-hand side of Fig. \ref{local_cycle} illustrates $\widetilde{H}_{1},\widetilde{H}_{2}$ and $\widetilde{H}_{3}$. We often describe such circumstances around $\widetilde{\gamma}$ (resp. $\widehat{\gamma}$) by thumbnailing each $\widetilde{H}_{i}$ (resp. $\widehat{H}_{i}$) with the ends as illustrated in Fig. \ref{local_cycle2}. We define the {\it interferency} of $\widehat{\gamma}$ as the number of all double points on $\widehat{\gamma}$ in $\widehat{G}$ which are not self double points of $\widehat{\gamma}$.
\section{Proof of Theorem \ref{main1}}

First we give two lemmas necessary for the proof of Theorem \ref{main1}.

\begin{Lemma}\label{interferency1}
Let $\varphi$ be a regular projection of a graph $G$, $f$ a spatial embedding of $G$ which projects on $\varphi$ and $\gamma$ a cycle of $G$. If $\varphi|_{\gamma}$ is trivial and the interferency of $\widehat{\gamma}$ is less than or equal to $1$, then there exists a trivialization disk for $f(\gamma)$.
\end{Lemma}

\begin{figure}[htbp]
      \begin{center}
\scalebox{0.4}{\includegraphics*{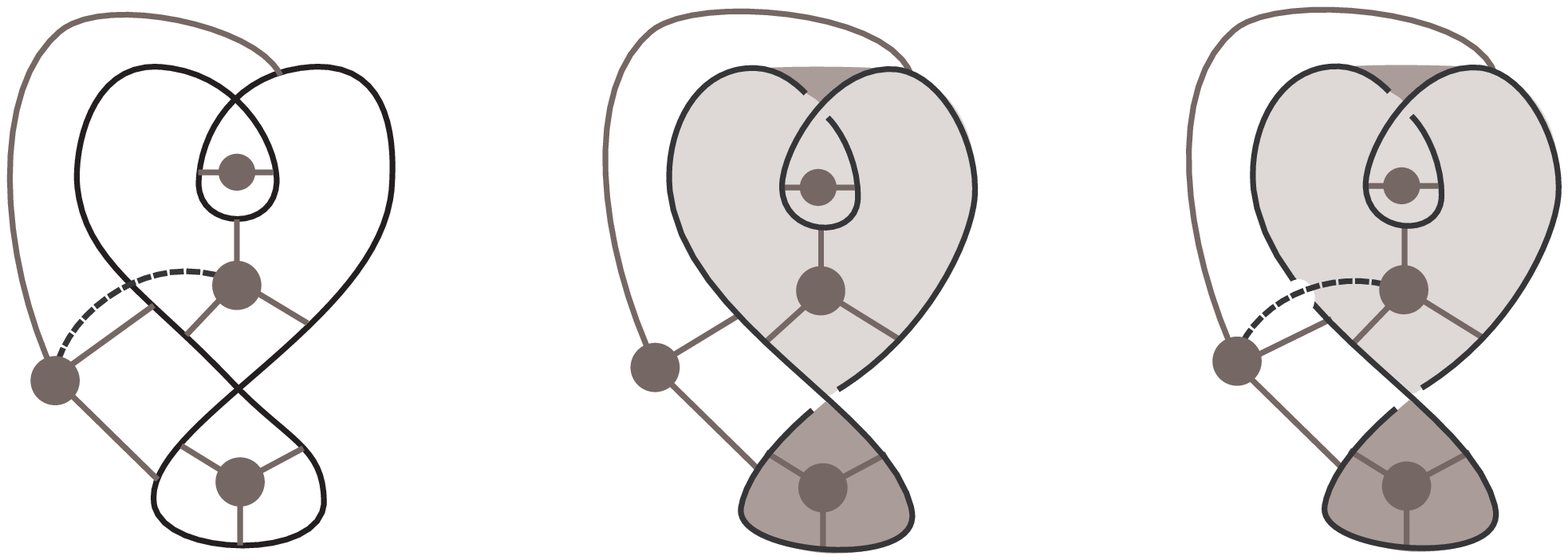}}
      \end{center}
   \caption{ }
  \label{trivial_disk}
\end{figure}

\begin{proof}
If the interferency of $\widehat{\gamma}$ is equal to $0$, we construct a canonical Seifert surface of $f(\gamma)$ from $\widetilde{\gamma}$ by applying the Seifert algorithm and, if necessary, isotope the surface so that it is located below each $f(H_i)$ with respect to the height defined by the natural projection $\pi$. Since $\varphi|_{\gamma}$ is trivial, the number of Seifert circles should be greater than the number of double points of $\varphi|_{\gamma}$ by one, which implies that the resulting surface is a disk. Therefore there exists a trivialization disk for $f(\gamma)$.

Now consider the case that the interferency of $\widehat{\gamma}$ is equal to $1$. If $\widehat{e}_{1}$ passes above ({\em resp.} under) $\widehat{\gamma}$, then we construct a canonical Seifert surface from $\widetilde{\gamma}$ so that it is located below ({\em resp.} above) each $f(H_i)$. Then we can obtain a trivialization disk for $f(\gamma)$, after isotpoing the Seifert surface (or $f(e_1)$ in relative sense) along the direction of the height so that $f(e_1)$ is above ({\em resp.} below) the resulting surface. Our construction is depicted in Fig. \ref{trivial_disk}.
\end{proof}

The following is a classification of regular projections of a cycle with at most three double points. See \cite[FIGURE 15]{Arnold96}.

\begin{figure}[htbp]
      \begin{center}
\scalebox{0.4}{\includegraphics*{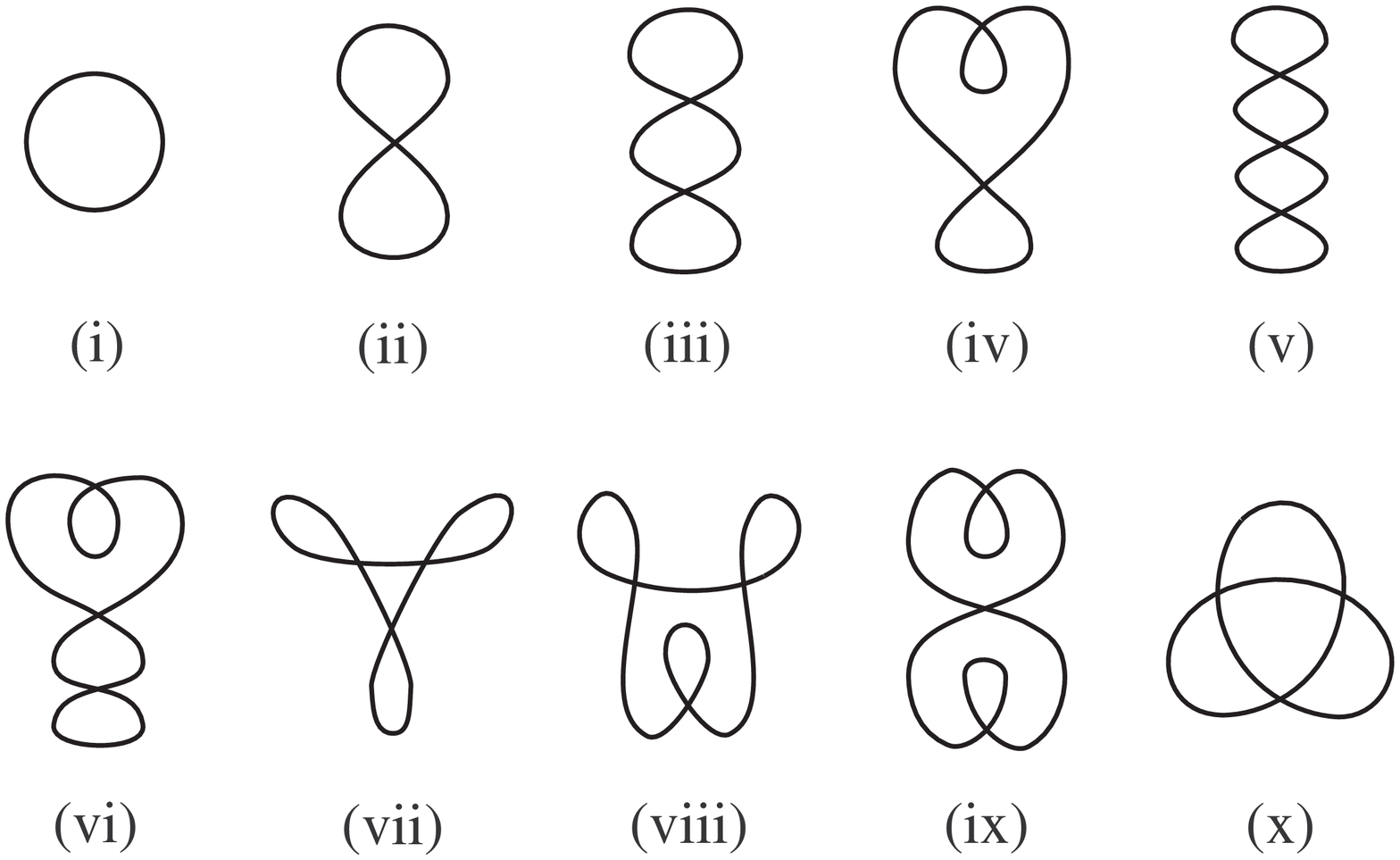}}
      \end{center}
   \caption{}
  \label{spherical_curves}
\end{figure}

\begin{Lemma}\label{id}
Let $\varphi$ be a regular projection of a cycle $\gamma$ with ${\rm cr}(\varphi)\le 3$. Then $\widehat{\gamma}$ is one of the ten projections as illustrated in Fig. \ref{spherical_curves} up to isotopy of ${\mathbb S}^{2}$.
\end{Lemma}

\noindent
\begin{proof}[Proof of Theorem \ref{main1}.]
Let $\gamma$ be a cycle of $G$. Let $e_{1},e_{2},\ldots,e_{m}\ (0\le m\le 3)$ be all different edges of $G$ which are not included in $\gamma$ such that $\widehat{\gamma}$ and $\widehat{e}_{i}$ produce double points of $\varphi$. By subdividing $G$ with some vertices of valency two if necessary, we may assume that $\widehat{\gamma}$ and $\widehat{e}_{i}$ produce exactly one double point of $\varphi$. We shall show that if $f$ does not contain a Hopf link or a trefoil knot then there exists a trivialization disk $D_{\gamma}$ for $f(\gamma)$. Then by Theorem \ref{keythm} we have the desired conclusion. Since ${\rm cr}(\varphi)\le 3$, $\widehat{\gamma}$ is one of the ten projections as illustrated in Fig. \ref{spherical_curves} up to isotopy of ${\mathbb S}^{2}$. Note that these regular projections are trivial except for (x). If $\widehat{\gamma}$ is any one of (iii), (iv), (v), (vi), (vii), (viii) or (ix), then the interferency of $\widehat{\gamma}$ is less than or equal to $1$ and by Lemma \ref{interferency1} there exists a trivialization disk for $f(\gamma)$. Thus, in the rest of the proof we show the claim in the case that $\widehat{\gamma}$ is (i), (ii) or (x).

Let us consider the case that $\widehat{\gamma}$ is (i) or (ii). If the interferency of $\widehat{\gamma}$ is less than or equal to $1$, then by Lemma \ref{interferency1} there exists a trivialization disk for $f(\gamma)$. So we assume that the interferency of $\widehat{\gamma}$ is $2$ or $3$. Since ${\rm cr}(\varphi)\le 3$, we may divide our situation about the circumstances around $\widehat{\gamma}$ into the four cases (1), (2), (3) and (4) as illustrated in Fig. \ref{disk2}. We remark here that there are ambiguities for positions of the ends of $\widehat{H}_{i}$ in Fig. \ref{disk2}, but they do not have an influence on our arguments except for the case (4d) as we will say later. In the following we observe a regular diagram of a spatial embedding of $G$ which projects on (1), (2), (3) or (4).
\begin{figure}[htbp]
      \begin{center}
\scalebox{0.45}{\includegraphics*{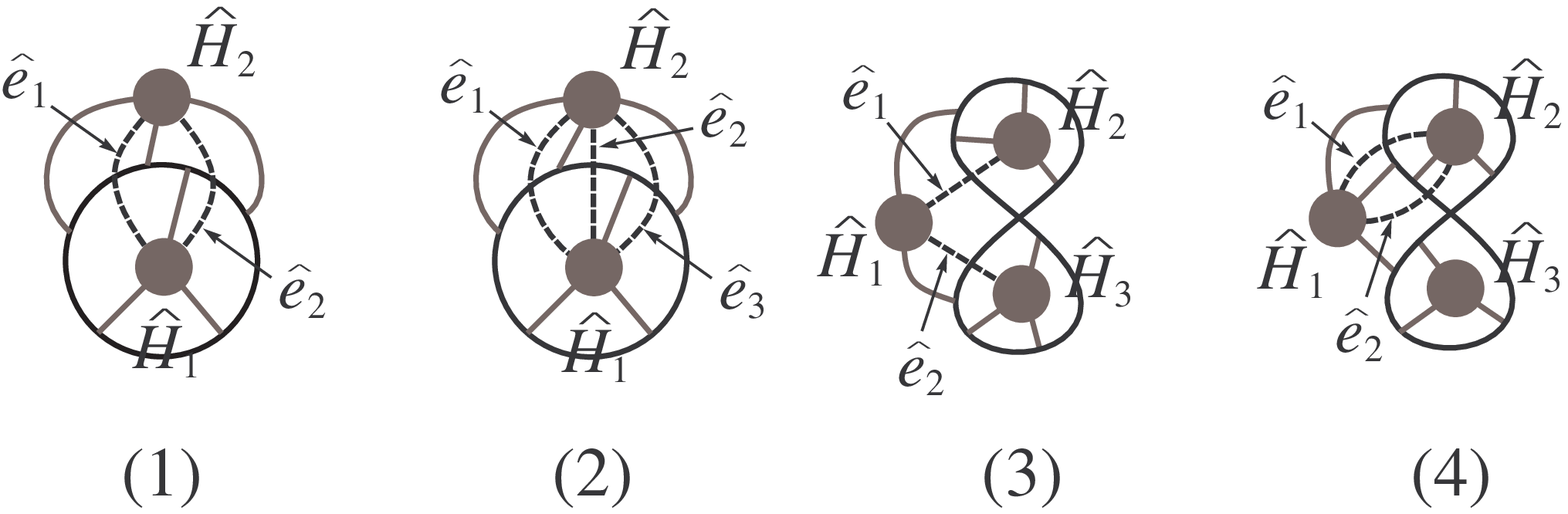}}
      \end{center}
   \caption{}
  \label{disk2}
\end{figure}

\noindent
(1) It is sufficient to consider the two cases (1a) and (1b) as illustrated in Fig. \ref{disk1}. The other cases can be shown by considering the mirror image embedding in the same way as the proof of Lemma \ref{interferency1} (after this we often adopt this argument and do not touch on it one by one). In the case (1a), it is clear that there exists a trivializing disk for $f(\gamma)$, see Fig. \ref{disk1}. Next we consider the case (1b). Since $f$ does not contain a Hopf link, we may assume that $e_{1}$ and $e_{2}$ each are incident to the different connected components of $H_{1}$ without loss of generality. Then we can see that there exists a trivializing disk for $f(\gamma)$, see Fig. \ref{disk1}.
\begin{figure}[htbp]
      \begin{center}
\scalebox{0.45}{\includegraphics*{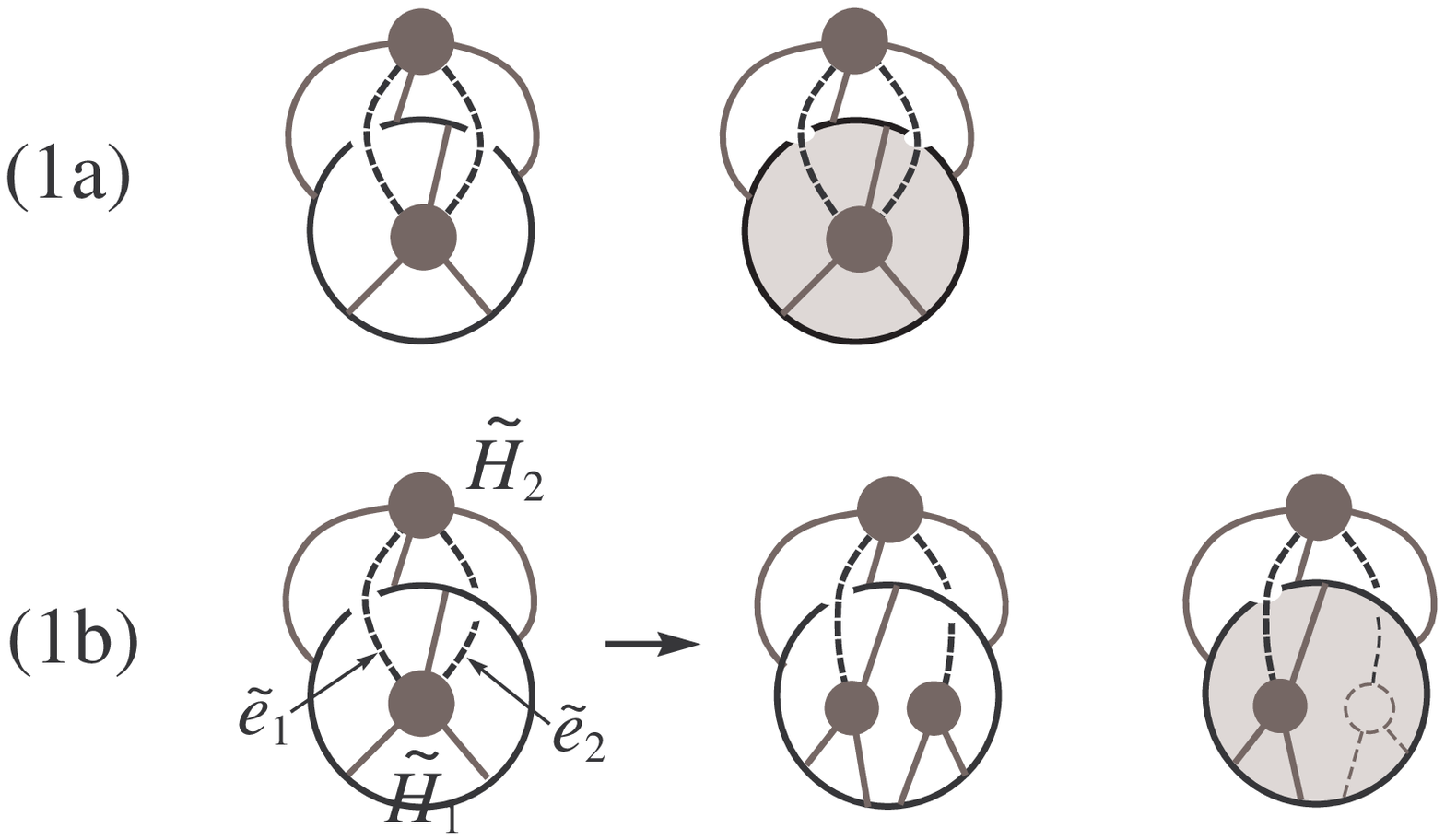}}
      \end{center}
   \caption{}
  \label{disk1}
\end{figure}

\noindent
(2) It is sufficient to consider the two cases (2a) and (2b) as illustrated in Fig. \ref{disk3}. In the case (2a), it is clear that there exists a trivializing disk for $f(\gamma)$, see Fig. \ref{disk3}. Next we consider the case (2b). Since $f$ does not contain a Hopf link, we may assume that both $e_{1}$ and $e_{2}$ are not incident to the connected component of $H_{1}$ to which $e_{3}$ is incident, or both $e_{1}$ and $e_{3}$ are not incident to the connected component of $H_{1}$ to which $e_{2}$ is incident without loss of generality. In the former case, it is clear that there exists a trivializing disk for $f(\gamma)$, see Fig. \ref{disk3}. In the latter case, we may assume that both $e_{1}$ and $e_{3}$ are incident to the same connected component of $H_{1}$. Since $f$ does not contain a Hopf link, we have that $e_{1}$ and $e_{3}$ are incident to the different connected components of $H_{2}$. If $e_{2}$ is not incident to the connected component of $H_{2}$ to which $e_{3}$ is incident, then we can see that there exists a trivializing disk for $f(\gamma)$, see Fig. \ref{disk3}. If $e_{2}$ is not incident to the connected component of $H_{2}$ to which $e_{1}$ is incident, then we also can see that there exists a trivializing disk for $f(\gamma)$, see Fig. \ref{disk3}.
\begin{figure}[htbp]
      \begin{center}
\scalebox{0.45}{\includegraphics*{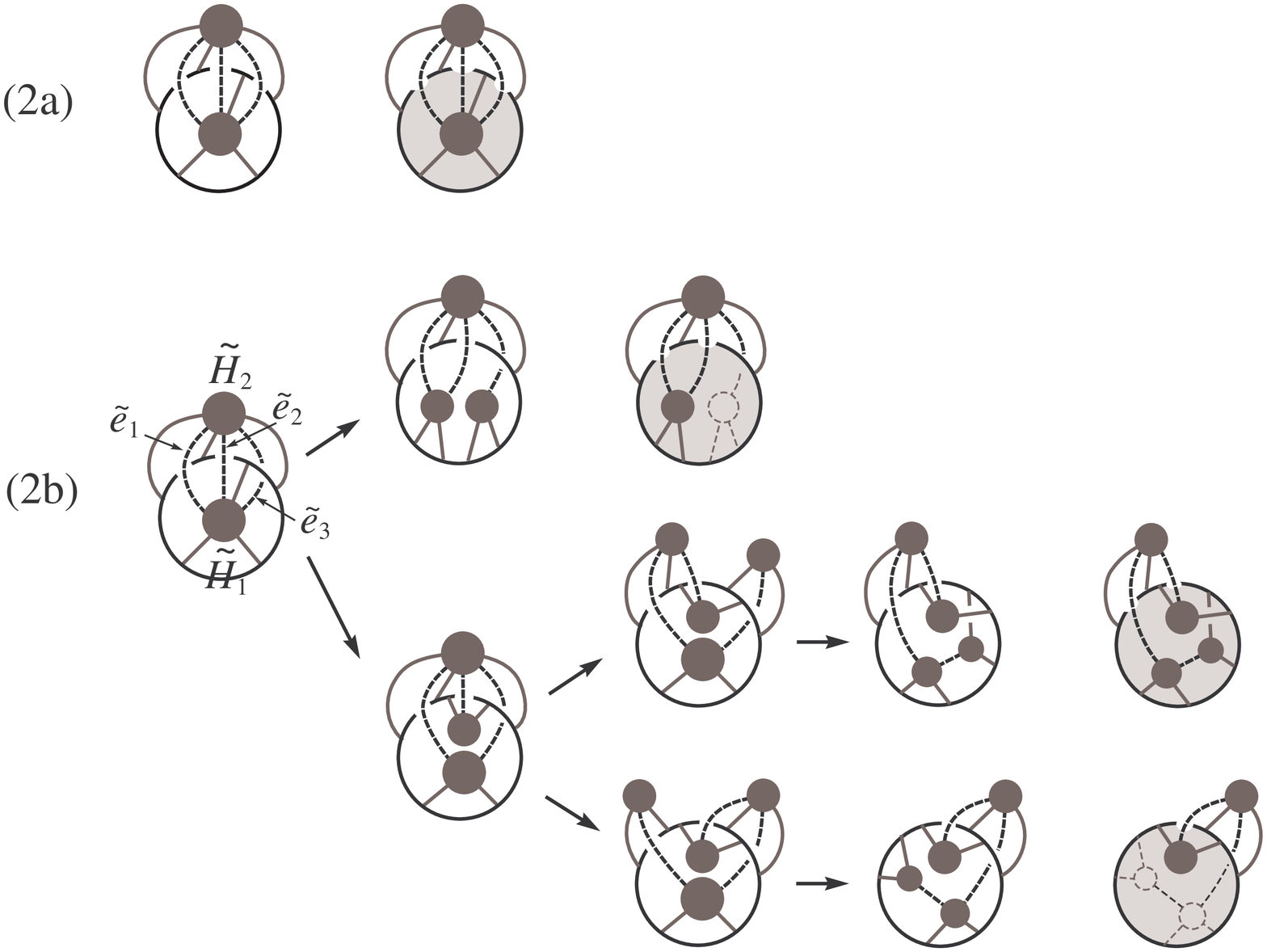}}
      \end{center}
   \caption{}
  \label{disk3}
\end{figure}

\noindent
(3) It is sufficient to consider the four cases (3a), (3b), (3c) and (3d) as illustrated in Fig. \ref{disk4}. In any cases we can see easily that there exists a trivializing disk for $f(\gamma)$, see Fig. \ref{disk4}.
\begin{figure}[htbp]
      \begin{center}
\scalebox{0.45}{\includegraphics*{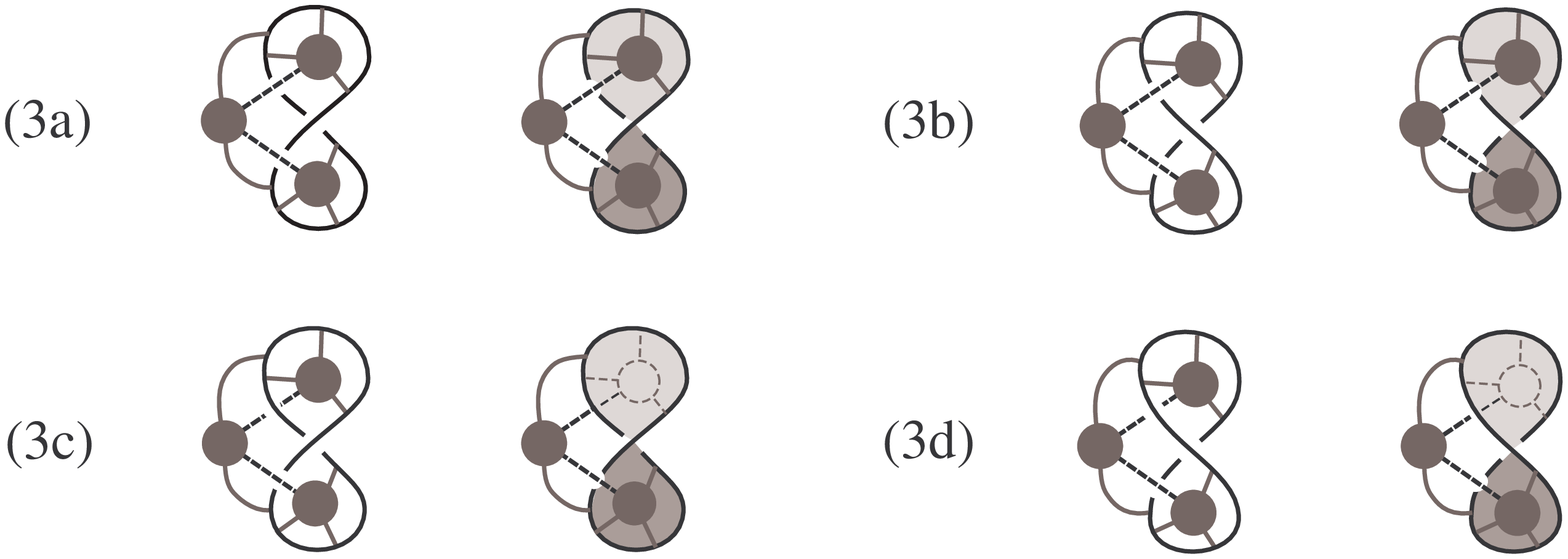}}
      \end{center}
   \caption{}
  \label{disk4}
\end{figure}

\noindent
(4) It is sufficient to consider the four cases (4a), (4b), (4c) and (4d) as illustrated in Fig. \ref{disk5}. In the cases (4a) and (4b), it is clear that there exists a trivializing disk for $f(\gamma)$, see Fig. \ref{disk6}. Next we consider the case (4c). Since $f$ does not contain a Hopf link, we have that $e_{1}$ and $e_{2}$ are incident to the different connected components of $H_{2}$, or $e_{1}$ and $e_{2}$ are incident to the different connected components of $H_{1}$. In either cases we can see that there exists a trivializing disk for $f(\gamma)$, see Fig. \ref{disk7}. Next we consider the case (4d). Since $f$ does not contain a Hopf link, we have that $e_{1}$ and $e_{2}$ are incident to the different connected components of $H_{2}$, or $e_{1}$ and $e_{2}$ are incident to the different connected components of $H_{1}$. In the former case, it is clear that there exists a trivializing disk for $f(\gamma)$, see Fig. \ref{disk8}. In the latter case, we may assume that both $e_{1}$ and $e_{2}$ are incident to the same connected component of $H_{2}$. We denote the connected component of $H_{1}$ to which $e_{i}$ is incident by $H_{1}^{(i)}\ (i=1,2)$. If there exist an end of $\widetilde{H}_{1}^{(1)}$ and an end of $\widetilde{H}_{1}^{(2)}$ each of which attaches to the boundary of $R_{3}$, then they must have a common vertex on the boundary of $R_{3}$ because $f$ does not contain a trefoil knot, see Fig. \ref{disk9}. Then we can see that there exists a trivializing disk for $f(\gamma)$, see Fig. \ref{disk8}.
\begin{figure}[htbp]
      \begin{center}
\scalebox{0.45}{\includegraphics*{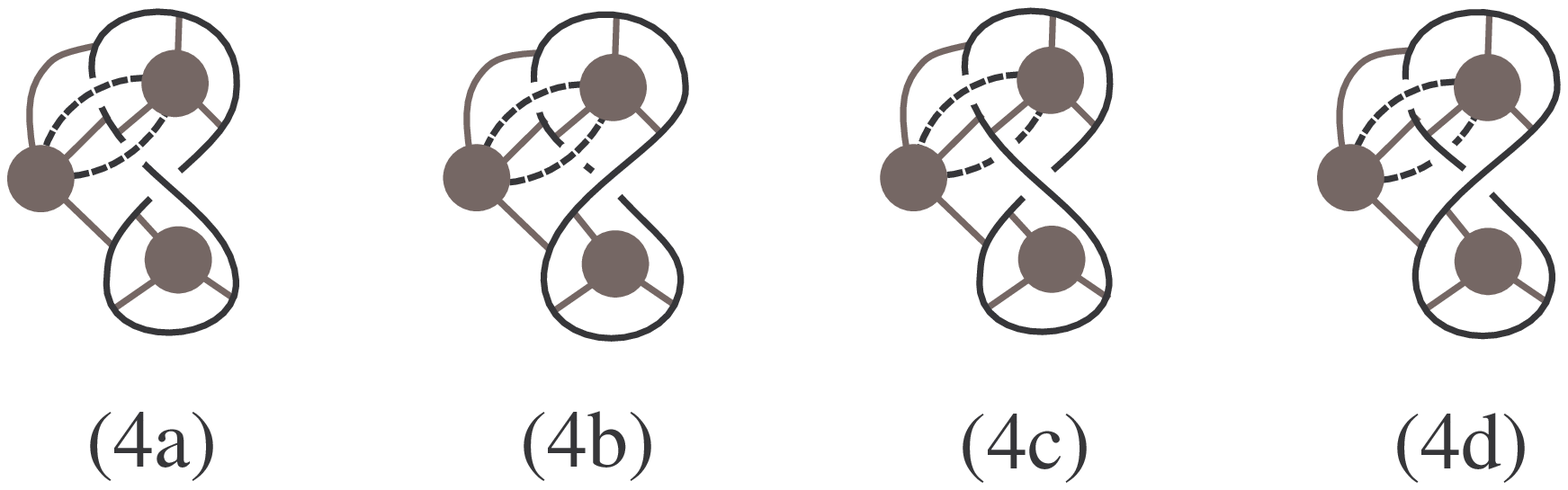}}
      \end{center}
   \caption{}
  \label{disk5}
\end{figure}
\begin{figure}[htbp]
      \begin{center}
\scalebox{0.45}{\includegraphics*{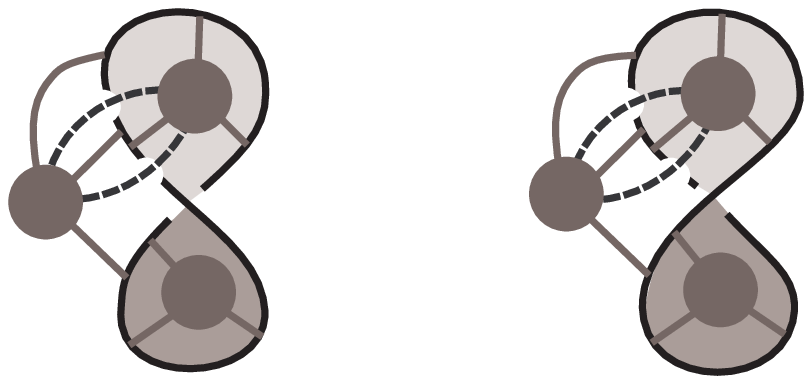}}
      \end{center}
   \caption{}
  \label{disk6}
\end{figure}
\begin{figure}[htbp]
      \begin{center}
\scalebox{0.45}{\includegraphics*{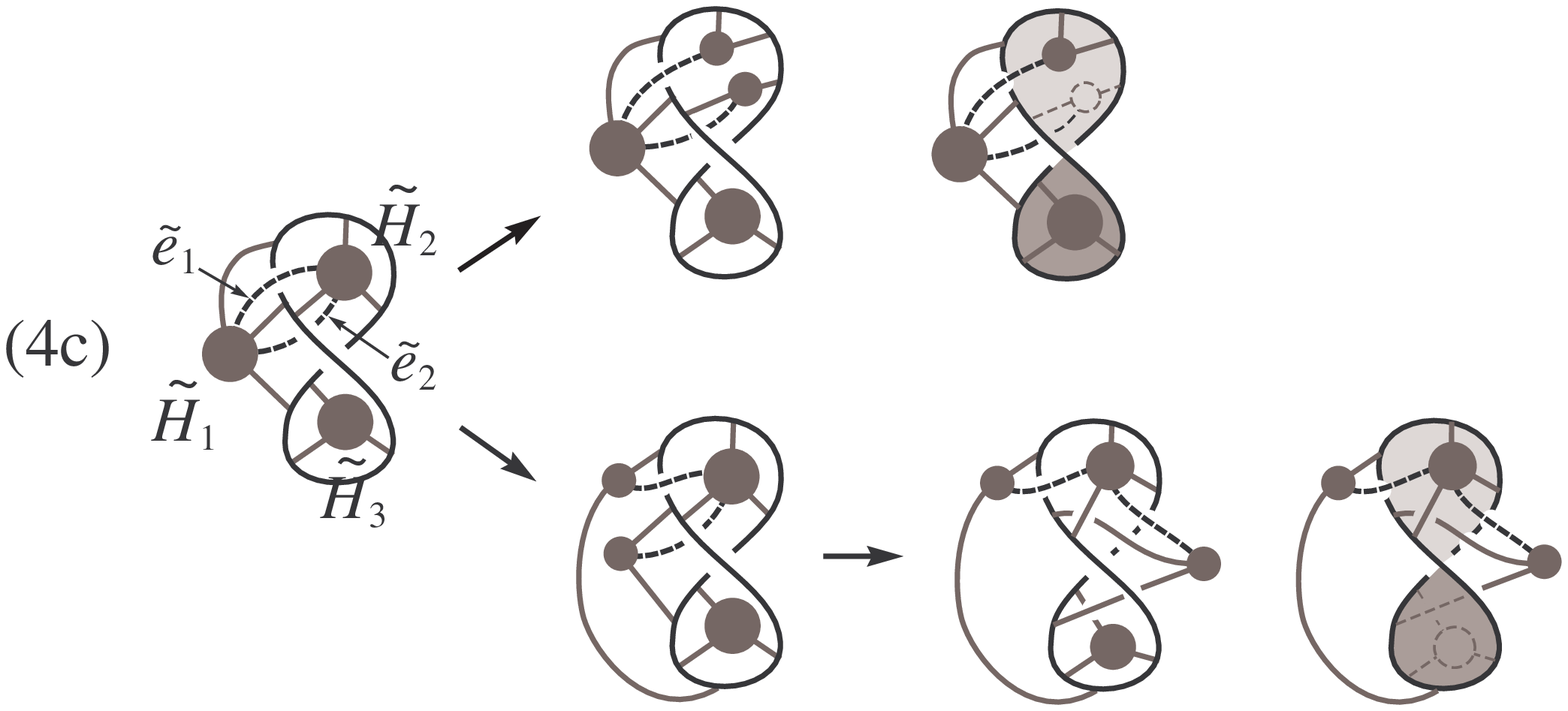}}
      \end{center}
   \caption{}
  \label{disk7}
\end{figure}
\begin{figure}[htbp]
      \begin{center}
\scalebox{0.44}{\includegraphics*{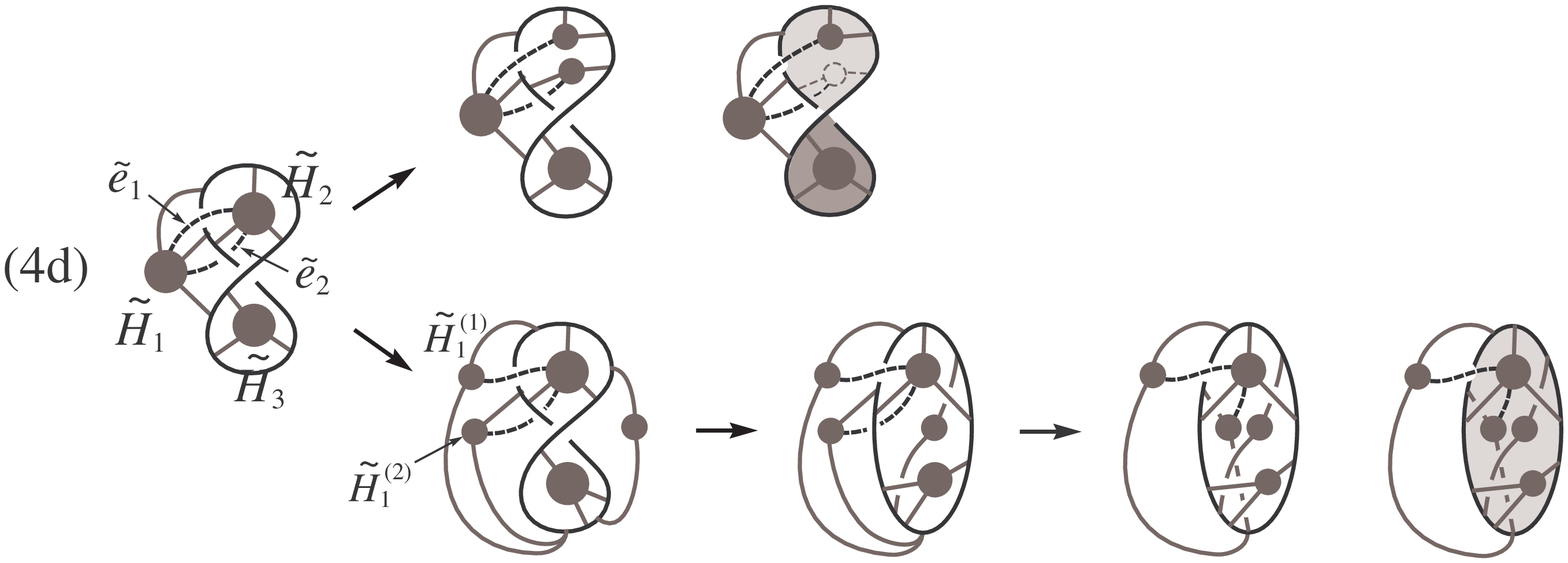}}
      \end{center}
   \caption{}
  \label{disk8}
\end{figure}
\begin{figure}[htbp]
      \begin{center}
\scalebox{0.45}{\includegraphics*{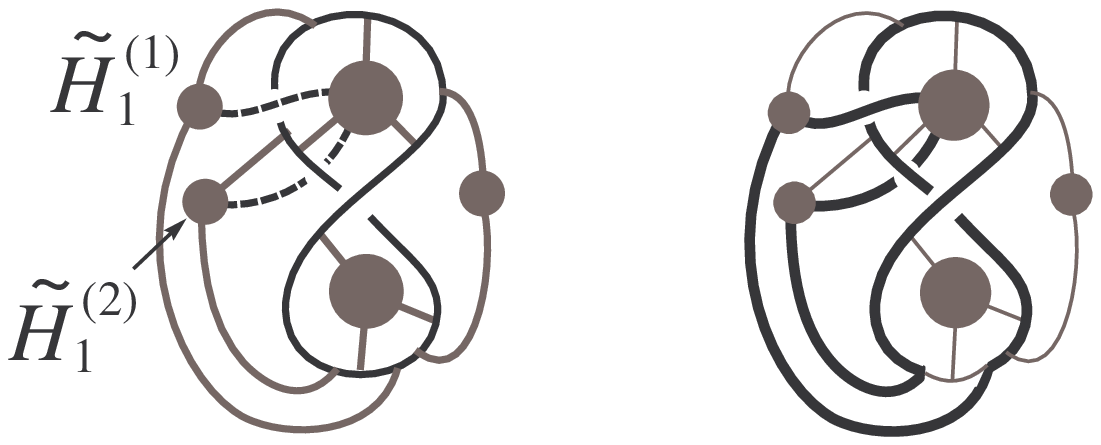}}
      \end{center}
   \caption{}
  \label{disk9}
\end{figure}

Finally let us consider the case that $\widehat{\gamma}$ is (x). Then the interferency of $\widehat{\gamma}$ is equal to $0$ and only trivial knots or trefoil knots project on it. Since $f$ does not contain a trefoil knot, we have that $f(\gamma)$ is a trivial knot. Thus our situation about the circumstances around $\widehat{\gamma}$ can be depicted as Fig. \ref{disk10}. And we can find a trivialization disk for $f(\gamma)$, see  Fig. \ref{disk10}. This completes the proof.
\begin{figure}[htbp]
      \begin{center}
\scalebox{0.45}{\includegraphics*{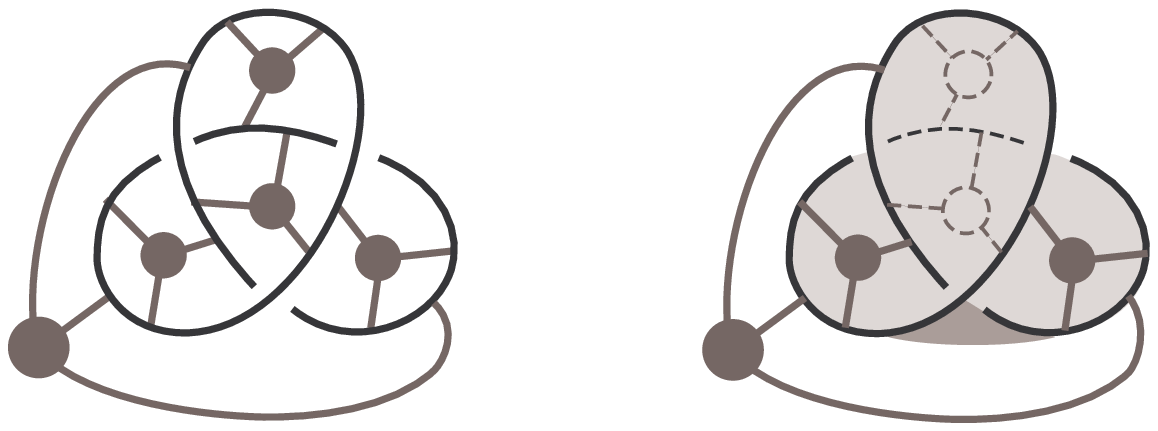}}
      \end{center}
   \caption{}
  \label{disk10}
\end{figure}
\end{proof}

\section{Proof of Theorem \ref{main2}}

In this section we prove Theorem \ref{main2}. For a regular projection $\varphi$ of a graph $G$, we denote the set of all equivalence classes of spatial embeddings of $G$ which project on $\varphi$ by ${\rm SE}(\varphi)$. We say that two regular projections $\varphi$ and $\psi$ of $G$ are {\it {\rm SE}-equivalent} ($\varphi \sim_{{\rm SE}} \psi$) if ${\rm SE}(\varphi)={\rm SE}(\psi)$.

\begin{proof}[Proof of Theorem \ref{main2}.]
(1) It is clear by Corollary \ref{main1cor} (1).

\noindent
(2) Let $\varphi$ be a regular projection of $G$ with ${\rm cr}(\varphi)=2$. If there exists a non-trivial spatial embedding $f$ of $G$ which projects on $\varphi$, then by Corollary \ref{main1cor} (2), $f$ contains a Hopf link. Since ${\rm cr}(\varphi)=2$, there exists a pair of disjoint cycles $\gamma$ and $\gamma'$ of $G$ such that $\widehat{\gamma}\cup \widehat{\gamma'}$ may be described as the left-hand side of Fig. \ref{Hopf_link0}. Then it is clear that each of the double points is Type-D. Therefore we have that if $\varphi$ has a double point of Type-S or Type-A then there does not exists a non-trivial spatial embedding of $G$ which projects on $\varphi$, namely $\varphi$ is trivial.
\begin{figure}[htbp]
      \begin{center}
\scalebox{0.45}{\includegraphics*{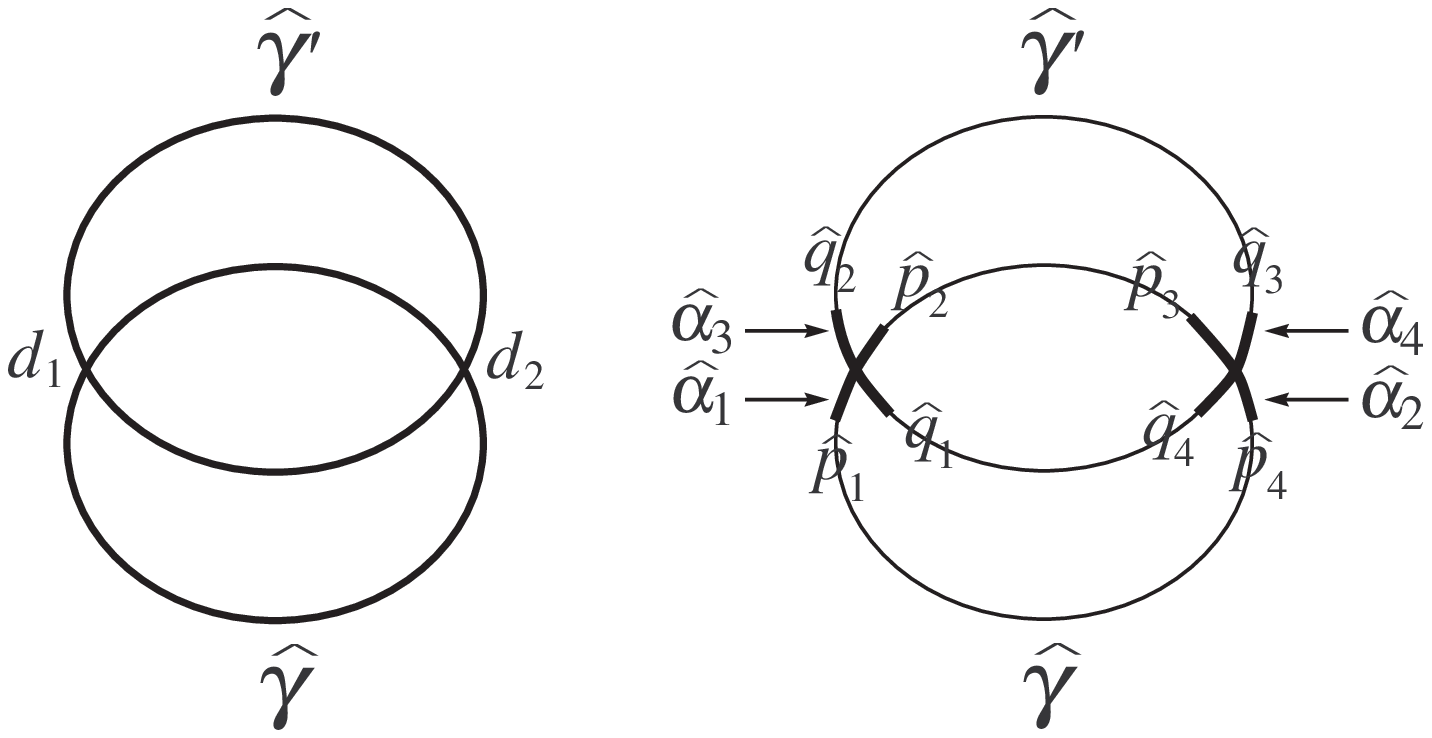}}
      \end{center}
   \caption{}
  \label{Hopf_link0}
\end{figure}

Now we assume that $\varphi$ is knotted. Let $d_{1}$ and $d_{2}$ be exactly two double points of $\varphi$. By considering sufficiently small compact neighborhoods $N_{1}$ and $N_{2}$ of $d_1$ and $d_2$ in ${\mathbb S}^2$, respectively, we can obtain disjoint simple subarcs $\alpha_1$, $\alpha_2$, $\alpha_3$ and $\alpha_4$ of $G$ each of which does not contain any vertex of $G$, so that $\varphi^{-1}(N_1) = \alpha_1 \cup \alpha_3$, $\varphi^{-1}(N_2) = \alpha_2 \cup \alpha_4$ and $\alpha_1 \cup \alpha_2 \subset \gamma$, $\alpha_3\cup \alpha_4 \subset \gamma'$. We put $\partial \alpha_{1}=\left\{p_{1},p_{2}\right\}$, $\partial \alpha_{2}=\left\{p_{3},p_{4}\right\}$, $\partial \alpha_{3}=\left\{q_{1},q_{2}\right\}$ and $\partial \alpha_{4}=\left\{q_{3},q_{4}\right\}$ so that $\widehat{p_{i}}$ and $\widehat{q}_{i}\ (i=1,2,3,4)$ are in the position as illustrated in the right-hand side of Fig. \ref{Hopf_link0}. We denote two arcs $\gamma\setminus({\rm int}\alpha_{1}\cup {\rm int}\alpha_{2})$ by $\iota_{1}$ and $\iota_{2}$ so that $\partial \iota_{1}=\left\{p_{2},p_{3}\right\}$ and $\partial \iota_{2}=\left\{p_{1},p_{4}\right\}$, and two arcs $\gamma'\setminus({\rm int}\alpha_{3}\cup {\rm int}\alpha_{4})$ by $\kappa_{1}$ and $\kappa_{2}$ so that $\partial \kappa_{1}=\left\{q_{2},q_{3}\right\}$ and $\partial \kappa_{2}=\left\{q_{1},q_{4}\right\}$. By giving over/under informations to $d_{1}$ and $d_{2}$ so that $\widetilde{\gamma}$ passes over $\widetilde{\gamma'}$, we can obtain the spatial embedding $f$ of $G$ which projects on $\varphi$ such that $f(\gamma\cup \gamma')$ is a trivial $2$-component link. Note that $f$ is non-trivial because $\varphi$ is knotted. Therefore by Corollary \ref{main1cor} (2), there exists a pair of disjoint cycles $\gamma''$ and $\gamma'''$ of $G$ such that $f(\gamma''\cup \gamma''')$ is a Hopf link. Since ${\rm cr}(\varphi)=2$, we may assume that $\alpha_{1}\cup \alpha_{4}\subset \gamma''$ and $\alpha_{2}\cup \alpha_{3}\subset \gamma'''$. Moreover, we may assume that there exists a pair of disjoint subarcs $\lambda_{1}$ and $\lambda_{2}$ of $\gamma''$ such that $\partial \lambda_{1}=\left\{p_{2},q_{3}\right\}$ and $\partial \lambda_{2}=\left\{p_{1},q_{4}\right\}$ without loss of generality. Then there exists a pair disjoint subarcs $\mu_{1}$ and $\mu_{2}$ of $\gamma'''$ such that $\partial \mu_{1}=\left\{p_{4},q_{2}\right\}$ and $\partial \mu_{2}=\left\{p_{3},q_{1}\right\}$. We denote the subgraph $\gamma\cup \gamma'\cup \gamma''\cup \gamma'''$ of $G$ by $H$. We call the closure of a connected component of $H\setminus {\gamma\cup \gamma'}$ in $H$ a {\it connector}. Note that a connector is a simple arc in $H$ whose boundary belongs to $\gamma\cup \gamma'$.

Since $H$ is planar, there exists an embedding $\psi:H\to {\mathbb S}^{2}$. For a subspace $S$ of $H$, we denote $\psi(S)$ also by $S$ as long as no confusion occurs. Then we may assume that $\gamma\cup \gamma'$ is positioned into ${\mathbb S}^{2}$ by $\psi$ as illustrated in Fig. \ref{shapes} (1) or (2). In any of the two cases, $\gamma'$ bounds a 2-disk $D'$ in ${\mathbb S}^{2}$ whose interior does not contain $\gamma$, and $\gamma$ bounds a $2$-disk $D$ in ${\mathbb S}^{2}$ whose interior contains $D'$. We denote the annulus $D\setminus {\rm int}D'$ by $A$ and the $2$-disk ${\mathbb S}^{2}\setminus {\rm int}D$ by $D''$. Note that there does not exist a connector between $\iota_{1}$ and $\iota_{2}$ (resp. $\kappa_{1}$ and $\kappa_{2}$) because if such a connector exists then ${\rm cr}(\varphi)> 2$. Therefore, if there exists a connector $c$ in $D''$ (resp. $D'$) then $\partial c\subset \iota_{i}$ (resp. $\partial c\subset \kappa_{i}$) for some $i$. Then we may assume that there does not exist any connector in $D'$ and $D''$ by making a detour through outermost connectors in $D'$ and $D''$ if necessary.
\begin{figure}[htbp]
      \begin{center}
\scalebox{0.35}{\includegraphics*{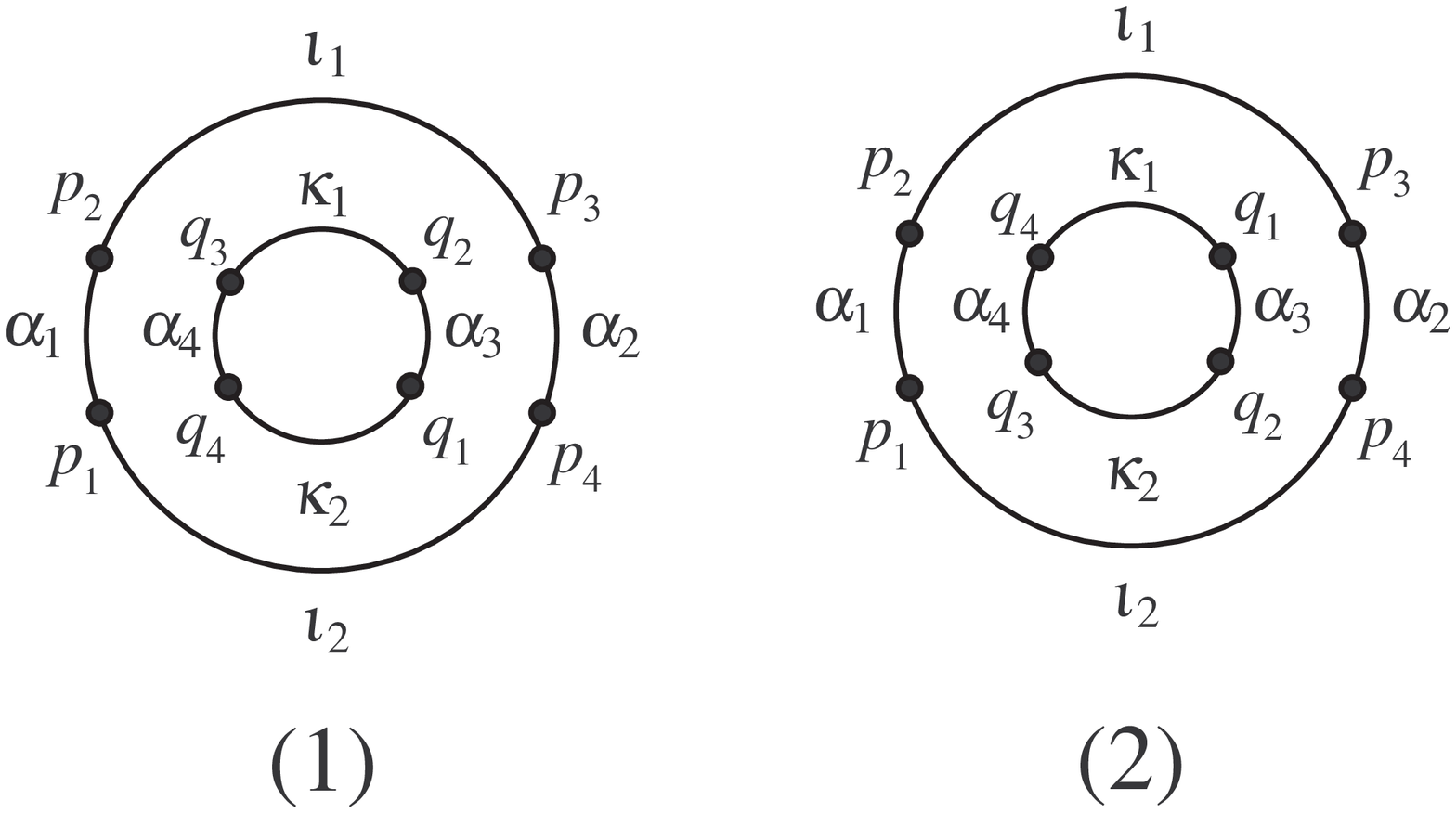}}
      \end{center}
   \caption{}
  \label{shapes}
\end{figure}

Now let us consider the case (1) of Fig. \ref{shapes}. Since there does not exist any connector in $D'$ and $D''$, we see that $\gamma''$ runs in $A$. Then we peel $\gamma''$ from $(\gamma\cup \gamma')\setminus (\alpha_{1}\cup \alpha_{4})$ in $A$ by applying local deformations as illustrated in Fig. \ref{peeling}. We also peel $\gamma'''$ from $(\gamma\cup \gamma')\setminus (\alpha_{2}\cup \alpha_{3})$ in the same way. By this operation, we can obtain new plane graph $H'$ from $H$. But it is easy to see that $H'$ contains a subspace which is homeomorphic to the complete bipartite graph on $3+3$ vertices, namely $H'$ is non-planar, see Fig. \ref{shapes2}. It is a contradiction. We can see that the case (2) of Fig. \ref{shapes} also yields a contradiction in a similar way. Hence we have that $\varphi$ is not knotted.
\begin{figure}[htbp]
      \begin{center}
\scalebox{0.4}{\includegraphics*{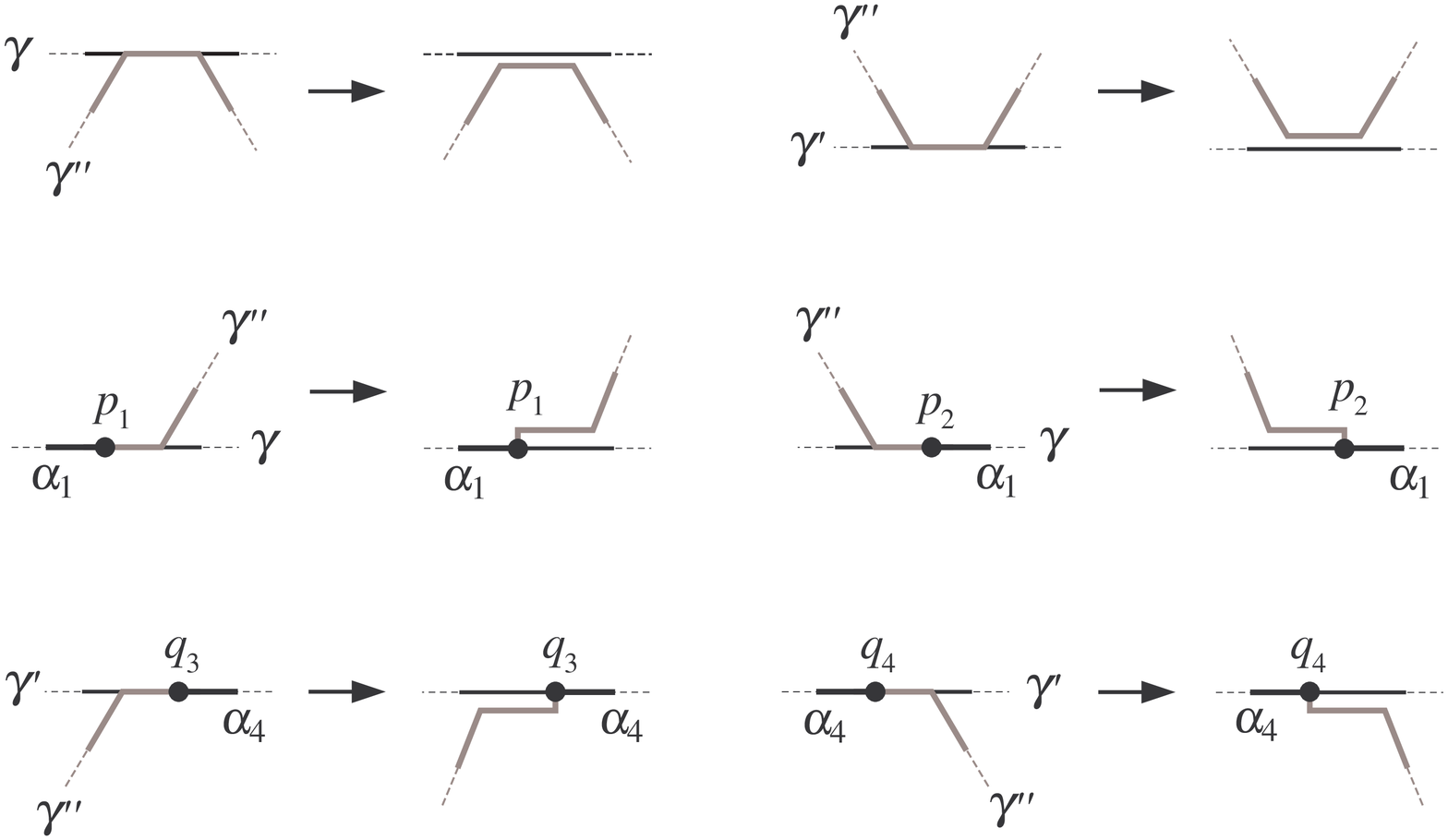}}
      \end{center}
   \caption{}
  \label{peeling}
\end{figure}
\begin{figure}[htbp]
      \begin{center}
\scalebox{0.35}{\includegraphics*{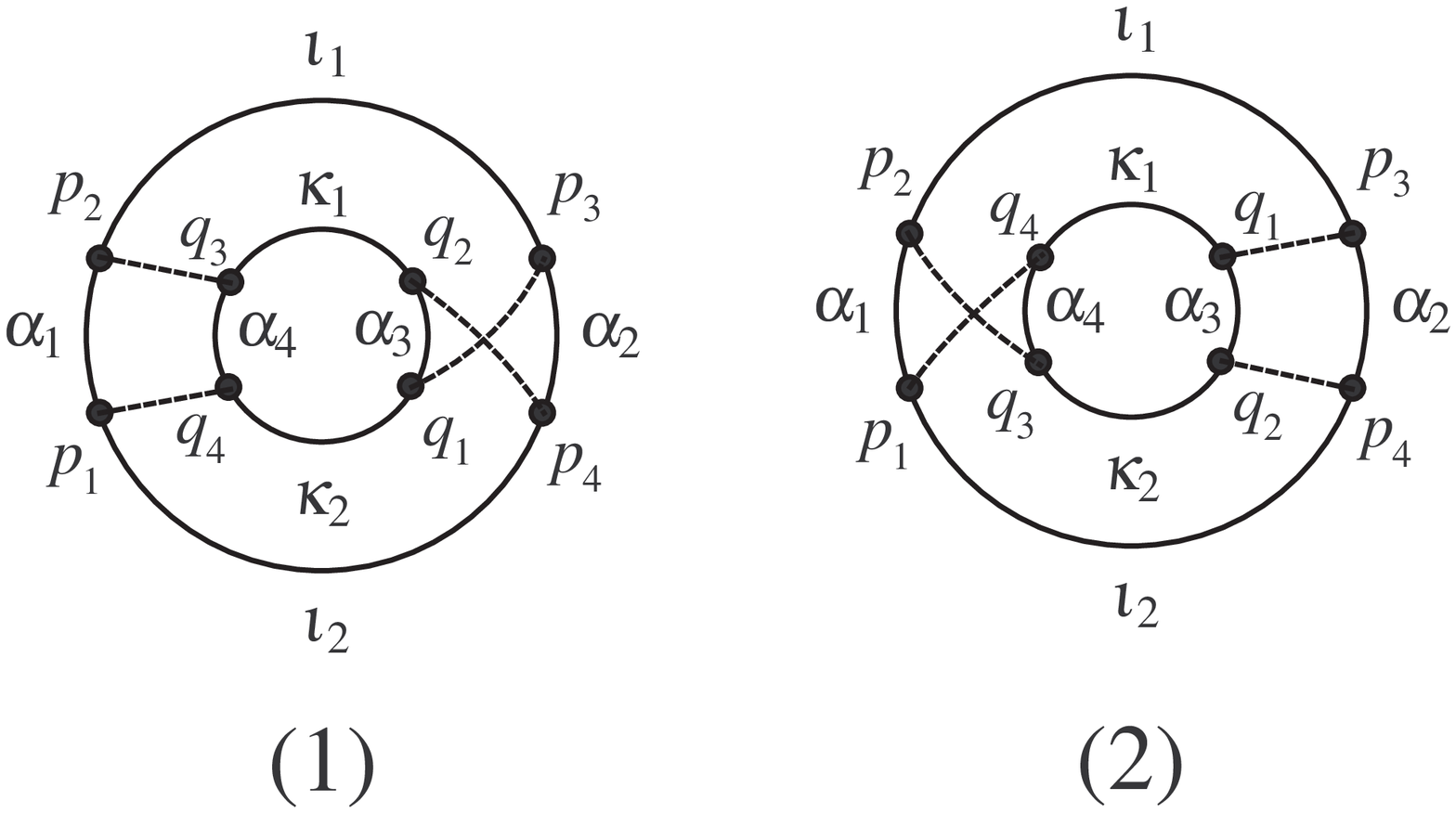}}
      \end{center}
   \caption{}
  \label{shapes2}
\end{figure}

\noindent
(3) Let $\varphi$ be a regular projection of $G$ with ${\rm cr}(\varphi)=3$. If $\varphi$ has a double point of Type-S, then we may divide our situation into the three cases (a), (b) and (c) as illustrated in Fig. \ref{not_knotted1}. In (a) and (b), $\varphi$ is {\rm SE}-equivalent to a regular projection $\psi$ of $G$ with ${\rm cr}(\psi)=2$ as illustrated in Fig. \ref{not_knotted1}. Note that $\psi$ is not knotted by Theorem \ref{main2} (2). Thus we have that $\varphi$ is not knotted. In (c), there exists a trivial spatial embedding of $G$ which projects on $\varphi$, see Fig. \ref{not_knotted1}. Thus we have that $\varphi$ is not knotted.
\begin{figure}[htbp]
      \begin{center}
\scalebox{0.4}{\includegraphics*{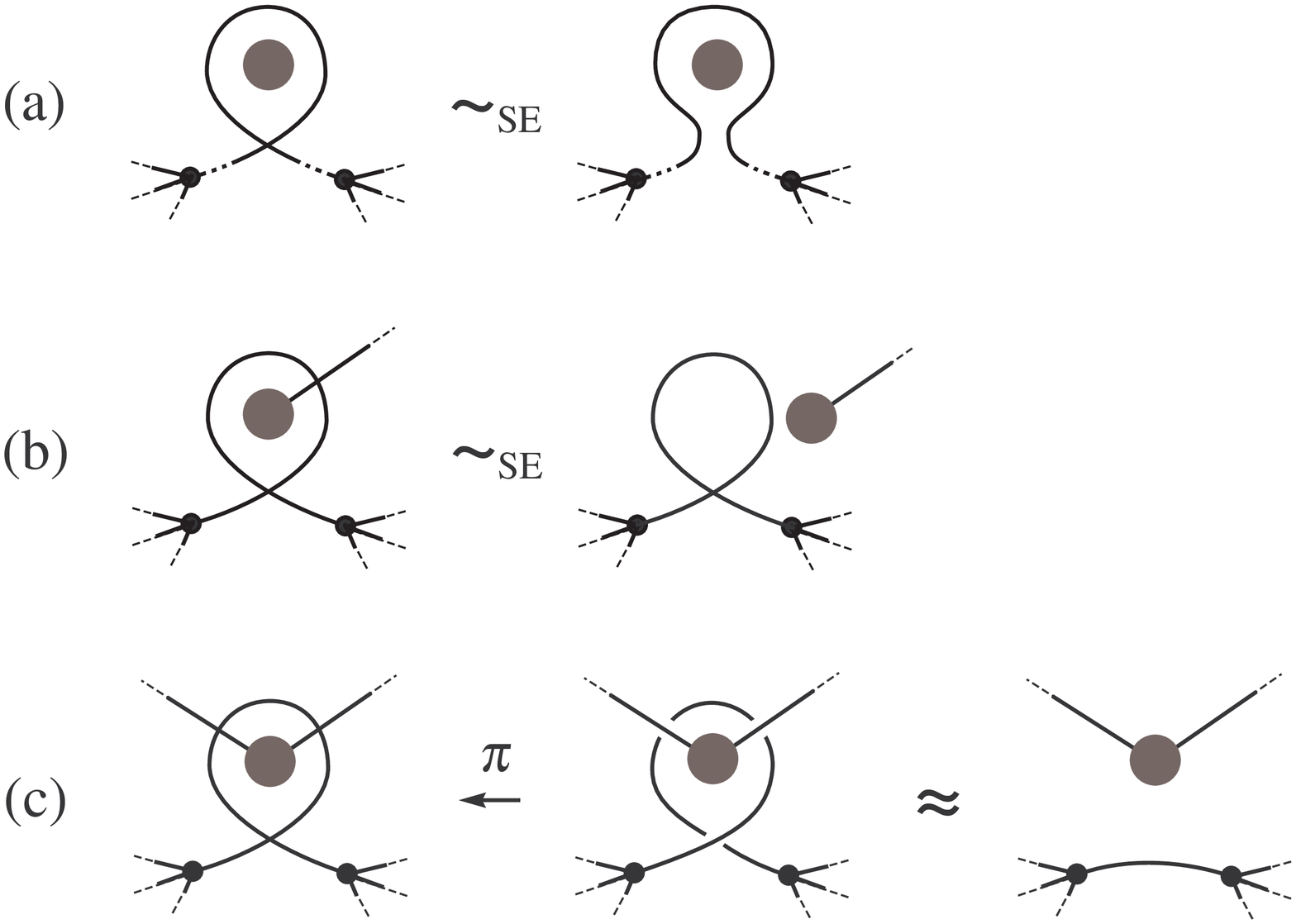}}
      \end{center}
   \caption{}
  \label{not_knotted1}
\end{figure}

If $\varphi$ does not have a double point of Type-S but has a double point of Type-A, then we may divide our situation into the four cases (e), (f), (g) and (h) as illustrated in Fig. \ref{not_knotted2}. In (e) and (f), we can show that $\varphi$ is not knotted in a similar way as (a) and (b). In (g), we can also show that $\varphi$ is not knotted in a similar way as (c). In (h), let $\psi$ be a regular projection of $G$ which is obtained from $\varphi$ by smoothing the double point of Type-A as illustrated in Fig. \ref{not_knotted3}. Since ${\rm cr}(\psi)=2$, we have that $\psi$ is not knotted by Theorem \ref{main2} (2). Note that any spatial embedding of $G$ which projects on $\psi$ also projects on $\varphi$, see Fig. \ref{not_knotted4}. Thus we have that if $\varphi$ is knotted then $\psi$ is also knotted. It is a contradiction. Hence we have that $\varphi$ is not knotted.
\begin{figure}[htbp]
      \begin{center}
\scalebox{0.4}{\includegraphics*{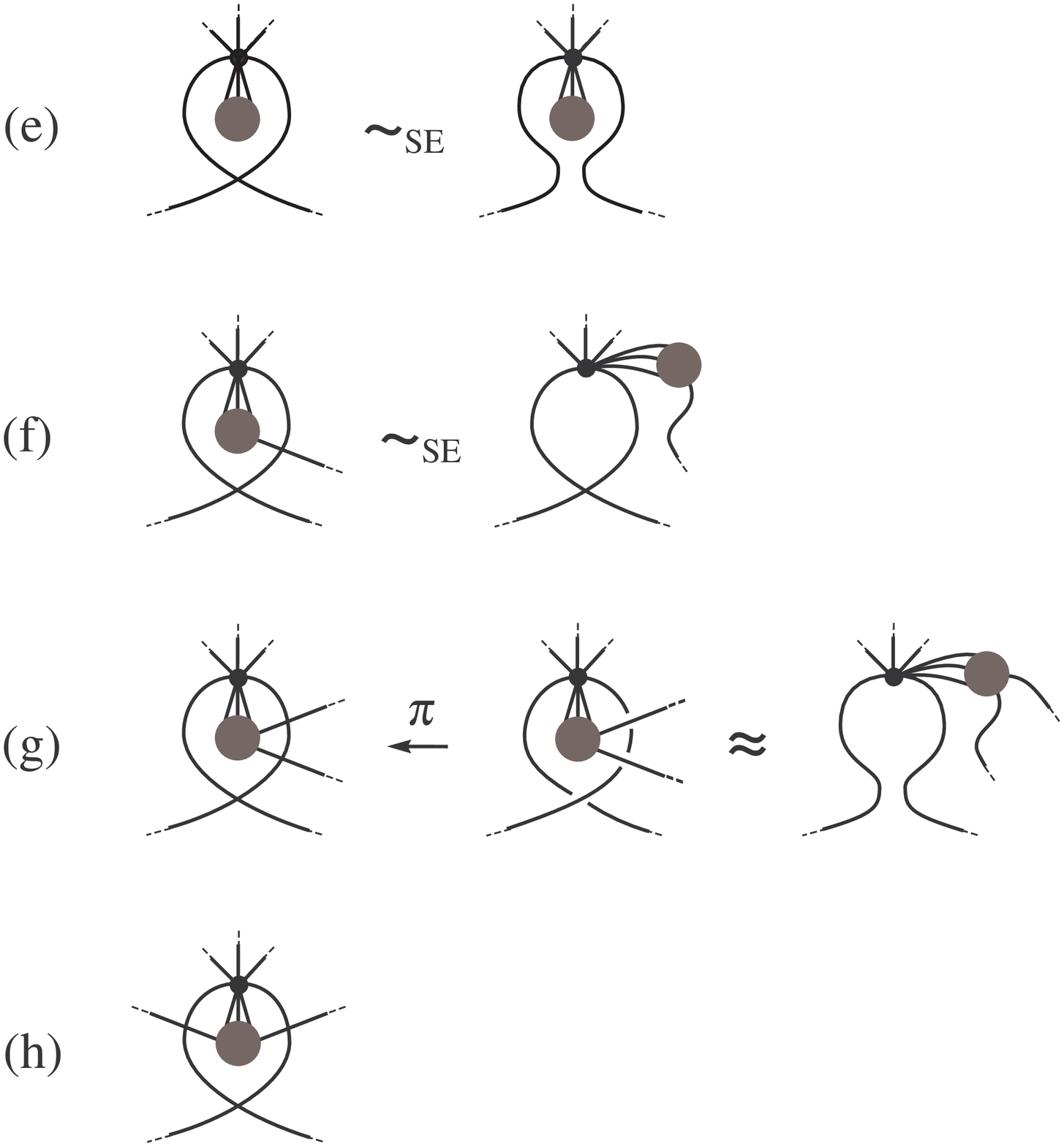}}
      \end{center}
   \caption{}
  \label{not_knotted2}
\end{figure}
\begin{figure}[htbp]
      \begin{center}
\scalebox{0.4}{\includegraphics*{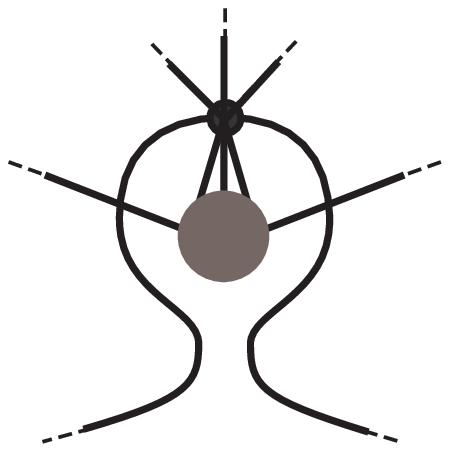}}
      \end{center}
   \caption{}
  \label{not_knotted3}
\end{figure}
\begin{figure}[htbp]
      \begin{center}
\scalebox{0.4}{\includegraphics*{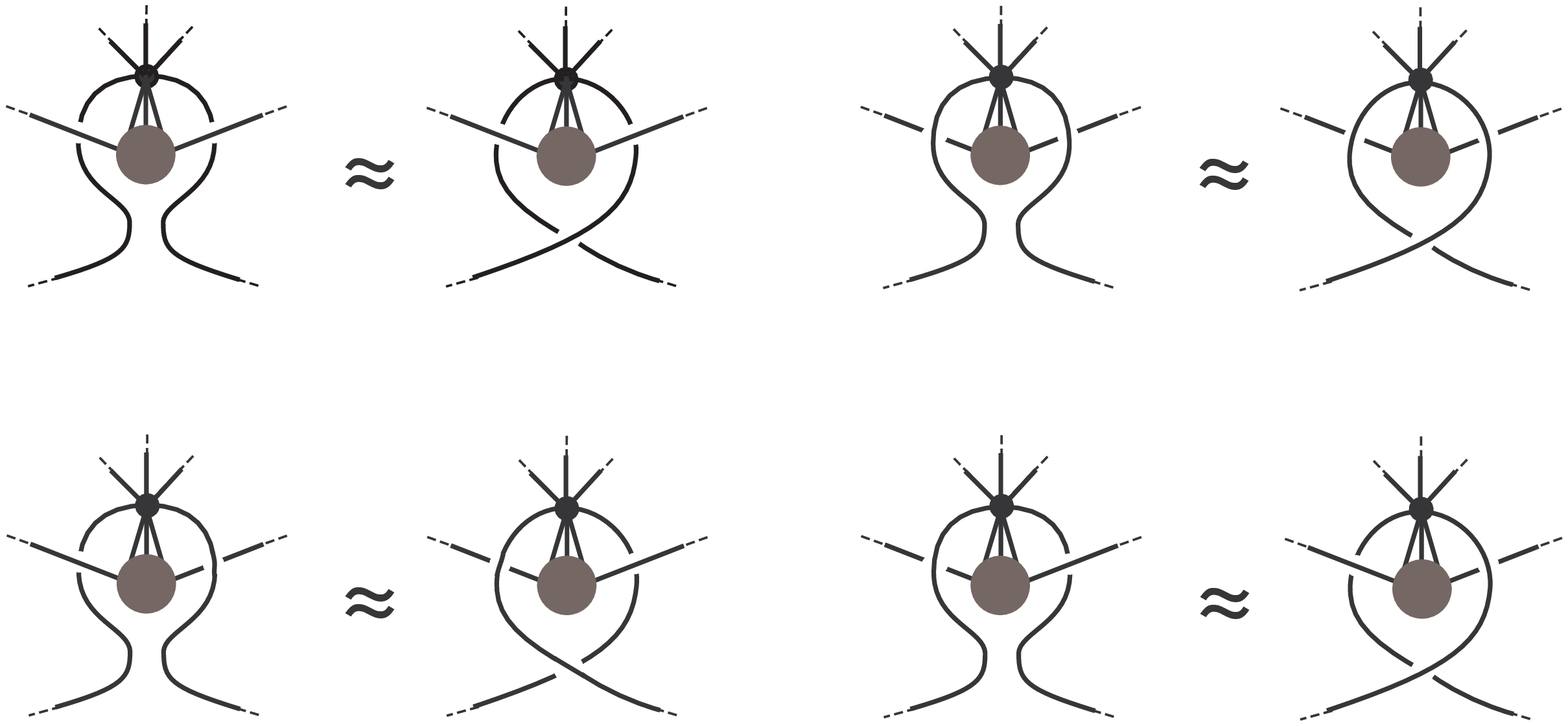}}
      \end{center}
   \caption{}
  \label{not_knotted4}
\end{figure}
\end{proof}

{\normalsize
}


\begin{thebibliography}{99}

\bibitem{Arnold96}
V. I. Arnold,
Remarks on the enumeration of plane curves,
{\it Topology of real algebraic varieties and related topics,} 17--32, Amer. Math. Soc. Transl. Ser. 2, {\bf 173}, {\it Amer. Math. Soc., Providence, RI,} 1996.

\bibitem{Kawauchi89}
A. Kawauchi,
Almost identical imitations of $(3,1)$-dimensional manifold pairs,
{\it Osaka J. Math.} {\bf 26} (1989), 743-758.

%\bibitem{Kura30} C. Kuratowski,
%Sur le probl\'eme des courbes gauches en topologie, (French)
%{\it Fund. Math.} {\bf 15} (1930), 271--283.

\bibitem{Lovasz06}
L. Lov{\'a}sz, 
Graph minor theory, 
{\it Bull. Amer. Math. Soc. (N.S.)} {\bf 43} (2006), 75--86 (electonic). 

%\bibitem{Mason69} W. K. Mason,
%Homeomorphic continuous curves in $2$-space are isotopic in $3$-space,
%{\it Trans. Amer. Math. Soc.} {\bf 142} (1969), 269--290.

\bibitem{NOTT05} R. Nikkuni, M. Ozawa, K. Taniyama and Y. Tsutsumi,
Newly found forbidden graphs for trivializability,
{\it J. Knot Theory Ramifications} {\bf 14} (2005), 523-538.

\bibitem{N07} R. Nikkuni,
Regular projections of spatial graphs,
{\it Knot Theory for Scientific Objects, Osaka City University Advanced Mathematical Institute Studies} {\bf 1} 111-128, {\it Osaka Municipal Universities Press}, 2007.

\bibitem{RS} N. Robertson and P. Seymour,
Graph minors XX. Wagner's conjecture, {\it J. Combin. Theory Ser. B} {\bf 92} (2004), 325--357.

%\bibitem{RST1}
%N. Robertson, P. Seymour and R. Thomas,
%Kuratowski chains,
%{\it J. Combin. Theory Ser. B} {\bf 64} (1995), 127-154.

%\bibitem{RST2}
%N. Robertson, P. Seymour and R. Thomas,
%Petersen family minors,
%{\it J. Combin. Theory Ser. B} {\bf 64} (1995), 155-184.

\bibitem{RST3}
N. Robertson, P. Seymour and R. Thomas,
Sachs' linkless embedding conjecture,
{\it J. Combin. Theory Ser. B} {\bf 64} (1995), 185-227.

\bibitem{ST91} M. Scharlemann and A. Thompson,
Detecting unknotted graphs in $3$-space,
{\it J. Diff. Geom.} {\bf 34} (1991), 539--560.

\bibitem{SS00} I. Sugiura and S. Suzuki,
On a class of trivializable graphs,
{\it Sci. Math.} {\bf 3} (2000), 193--200.

\bibitem{Tamu04} N. Tamura,
On an extension of trivializable graphs,
{\it J. Knot Theory Ramifications} {\bf 13} (2004), 211--218.

\bibitem{Tani95} K. Taniyama,
Knotted projections of planar graphs,
{\it Proc. Amer. Math. Soc.} {\bf 123} (1995), 3575--3579.

\bibitem{TT96} K. Taniyama and T. Tsukamoto,
Knot-inevitable projections of planar graphs,
{\it J. Knot Theory Ramifications} {\bf 5} (1996), 877-883.

%\bibitem{Tani-Yoshi98}
%K. Taniyama and C. Yoshioka,
%Regular projections of knotted handcuff graphs,
%{\it J. Knot Theory Ramifications} {\bf 7} (1998), 509--517.

\bibitem{Wu92}
Y. Q. Wu,
On planarity of graphs in $3$-manifolds,
{\it Comment. Math. Helv.} {\bf 67} (1992), 635--647.

\bibitem{Wu93}
Y. Q. Wu,
On minimally knotted embedding of graphs,
{\it Math. Z.} {\bf 214} (1993), 653--658.

\end{thebibliography}
\end{document}